\documentclass[10pt,a4paper]{amsart}

\usepackage[utf8x]{inputenc}
\usepackage{ucs}
\usepackage{tikz}
\usepackage{amsmath}
\usepackage{amsthm}
\usepackage{amsfonts}
\usepackage{amssymb}
\usepackage{hyperref}
\usepackage{dsfont}
\usepackage{mathrsfs}
\usepackage{eqnarray}
\write18{}
\usepackage{xcolor}
\usepackage{graphicx}
\usepackage{pdfpages} 
\usepackage[nothm]{thmbox}
\usepackage{esint}
\usepackage{comment}

\hypersetup{
    colorlinks=true,
    linkcolor=red,
    filecolor=magenta,      
    urlcolor=cyan,
}

\numberwithin{equation}{section}

\newcommand{\EE}{\mathbb{E}}

\newcommand{\PP}{\mathbb{P}}

\newcommand{\ZZ}{\mathbb{Z}}
\newcommand{\RR}{\mathbb{R}}

\newcommand{\NN}{\mathbb{N}}
\newcommand{\inte}[1]{\underset{#1}{\int}}

\newcommand{\infi}[1]{\underset{#1}{\inf}}
\newcommand{\supr}[1]{\underset{#1}{\sup}}
\newcommand{\somme}[1]{\underset{#1}{\sum}}

\newcommand{\tends}[1]{\underset{#1}{\longrightarrow}}

\DeclareMathOperator{\pmass}{\mathrm{pmass}}
\DeclareMathOperator{\mass}{\mathrm{mass}}

\DeclareMathOperator{\vol}{vol}

\DeclareMathOperator{\Id}{Id}

\DeclareMathOperator{\DR}{DR}
\DeclareMathOperator{\area}{Area}

\theoremstyle{definition}

\newtheorem{definition}[equation]{Definition}
\newtheorem{remark}[equation]{Remark}

\theoremstyle{theorem}

\newtheorem{theoreme}[equation]{Theorem}
\newtheorem{proposition}[equation]{Proposition}
\newtheorem{proposition*}[equation]{Proposition}

\newtheorem{lemma}[equation]{Lemma}

\title{A Cheeger-like inequality for coexact 1-forms}

\begin{document}

\begin{abstract}
	 We state and prove a Cheeger-like inequality for coexact 1-forms on closed orientable Riemannian manifolds.
\end{abstract}

\maketitle

\tableofcontents

\section{Introduction}
\label{sec.introduction}

Let $M$ be a closed connected Riemannian manifold of dimension $d$ and $\Delta = -\mbox{div}\circ \mbox{grad}$ its Laplacian, acting on $\mathcal{C}^\infty (M)$ (we will not explicitly refer to the metric tensor itself unless explicitly mentioned). The celebrated Cheeger's inequality relates the first non zero eigenvalue $\lambda _{1} (M)$ of the Laplacian 
with the 'Cheeger constant' $h(M)$ defined as 
$$ h(M) := \inf_{ S \hookrightarrow M, \ [S] = 0} \ \frac{\mu_{d-1}(S)}{\min(\mu(V_1), \mu(V_2))},$$
where the infinum is taken over all separating embeddings of codimension 1 submanifold $S$ into $M$, where $\mu_{d-1}$ [resp. $\mu$] stands for the $(d-1)$-dimensional Riemannian measure,
[resp. the Riemannian volume] of $M$ and where $V_1$ and $V_2$ are the two resulting connected components of $M \setminus S$. Cheeger's theorem states
\begin{theoreme}[Cheeger's Inequality]
	Let $M$ be a closed Riemannian manifold and let $\lambda_{1}(M)$ be the smallest eigenvalue of the Laplace operator acting on functions of $M$. Then 
	$$ \lambda_1 (M)  \ge \frac{h(M)^2}{4} \ .$$
\end{theoreme}
We can see the Cheeger constant as measuring a geometric cost of disconnecting the manifold $M$ and rephrase the above Theorem saying that if a closed manifold $M$ has a small non zero eigenvalue, then it tends to become disconnected. The heuristic behind such a statement is that a small eigenvalue of the Laplace operator may be thought of as an almost other dimension of $\ker \Delta$, which, by the maximum principle, has the same dimension as the numbers of connected components of $M$.

\begin{figure}[h!]
\begin{center}
	\def\svgwidth{0.8 \columnwidth}
%% Creator: Inkscape 1.0.1 (c497b03c, 2020-09-10), www.inkscape.org
%% Accompanies image file 'cheeger_alter.pdf' (pdf, eps, ps)
\begingroup%
  \makeatletter%
  \providecommand\color[2][]{%
    \errmessage{(Inkscape) Color is used for the text in Inkscape, but the package 'color.sty' is not loaded}%
    \renewcommand\color[2][]{}%
  }%
  \providecommand\transparent[1]{%
    \errmessage{(Inkscape) Transparency is used (non-zero) for the text in Inkscape, but the package 'transparent.sty' is not loaded}%
    \renewcommand\transparent[1]{}%
  }%
  \providecommand\rotatebox[2]{#2}%
  \newcommand*\fsize{\dimexpr\f@size pt\relax}%
  \newcommand*\lineheight[1]{\fontsize{\fsize}{#1\fsize}\selectfont}%
  \ifx\svgwidth\undefined%
    \setlength{\unitlength}{678.13705898bp}%
    \ifx\svgscale\undefined%
      \relax%
    \else%
      \setlength{\unitlength}{\unitlength * \real{\svgscale}}%
    \fi%
  \else%
    \setlength{\unitlength}{\svgwidth}%
  \fi%
  \global\let\svgwidth\undefined%
  \global\let\svgscale\undefined%
  \makeatother%
  \begin{picture}(1,0.4384026)%
    \lineheight{1}%
    \setlength\tabcolsep{0pt}%
    \put(0,0){\includegraphics[width=\unitlength,page=1]{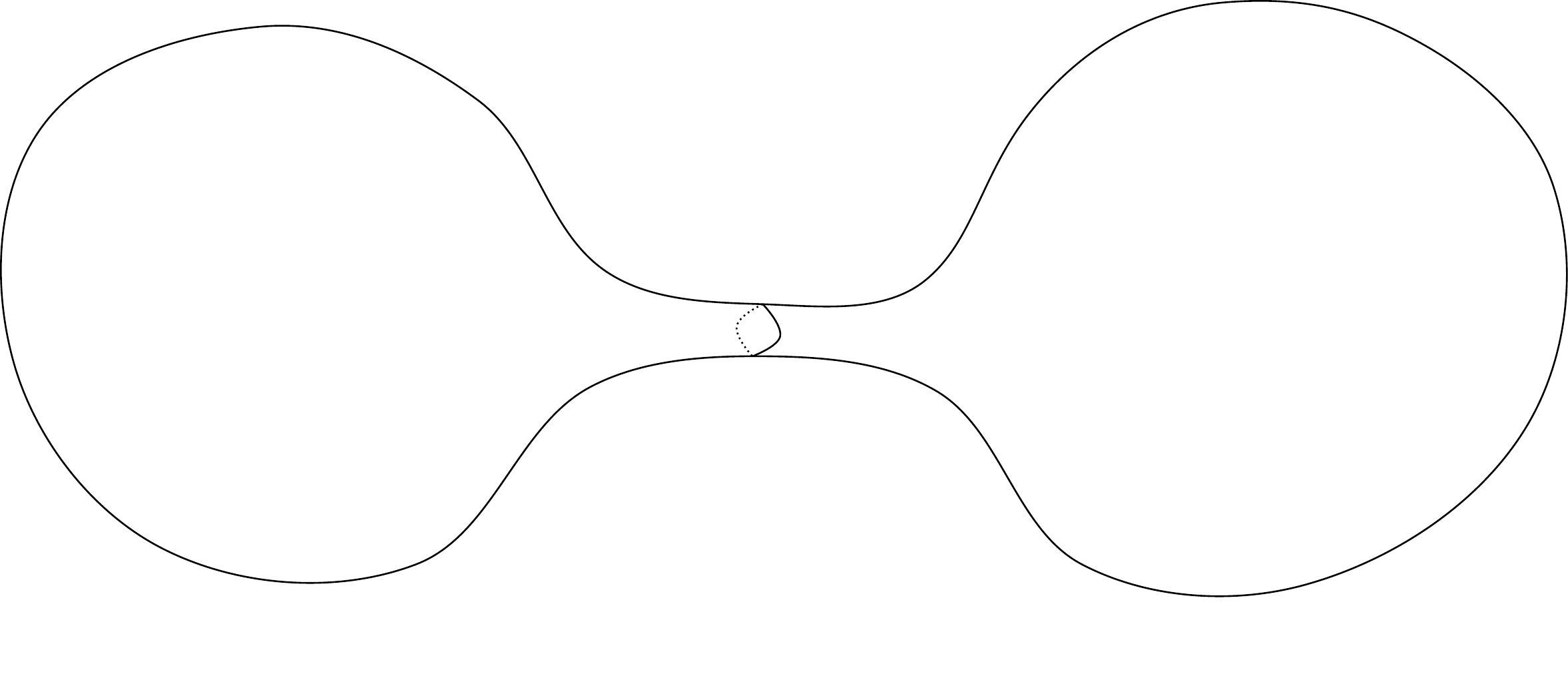}}%
    \put(0.10863247,0.22589455){\color[rgb]{0,0,0}\makebox(0,0)[lt]{\lineheight{1.25}\smash{\begin{tabular}[t]{l}$V_1$\end{tabular}}}}%
    \put(0.74998437,0.20556344){\color[rgb]{0,0,0}\makebox(0,0)[lt]{\lineheight{1.25}\smash{\begin{tabular}[t]{l}$V_2$\end{tabular}}}}%
    \put(0.44317051,0.15750832){\color[rgb]{0,0,0}\makebox(0,0)[lt]{\lineheight{1.25}\smash{\begin{tabular}[t]{l}$S$\end{tabular}}}}%
    \put(0.08090835,0.01334277){\color[rgb]{0.69803922,0,0}\makebox(0,0)[lt]{\lineheight{1.25}\smash{\begin{tabular}[t]{l}$f=-1$\end{tabular}}}}%
    \put(0.7666189,0.00779794){\color[rgb]{0.71764706,0,0}\makebox(0,0)[lt]{\lineheight{1.25}\smash{\begin{tabular}[t]{l}$f=1$\end{tabular}}}}%
    \put(0.41359811,0.27210145){\color[rgb]{0.76862745,0,0}\makebox(0,0)[lt]{\lineheight{1.25}\smash{\begin{tabular}[t]{l}$\nabla f \neq 0$\end{tabular}}}}%
  \end{picture}%
\endgroup%

	\end{center}
\caption{The 'dumb-bell' manifold represented here is almost disconnected. It carries a small eigenvalue: the function $f$ which values $+1$ on one side, $-1$ on the other side and whose gradient variations are concentrated around the hypersurface realising the infimum defining Cheeger's constant can be shown to have small energy.}
	\label{fig.cheeger}
\end{figure}
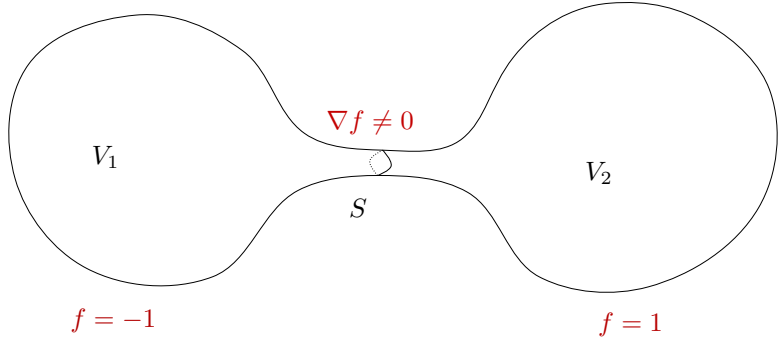

\bigskip
This article is dedicated to state and prove a Cheeger-like inequality for the first positive eigenvalue of the Hodge-Laplace operator acting on $1$-forms. The question of a Cheeger-like inequality for differential forms was raised in Cheeger's original article \cite{artcheegerspectre} and is formulated as Yau's problem 79 \cite{livreyauopenquestions}. The first positive eigenvalue of the Hodge Laplacian was relativity poorly studied compared to the functions case. Most of the work dedicated to the subject so far was turned toward the study of classes of examples, as constructions for which it exhibits a desired behaviour \cite{articlecolboiscourtoisspectrum1form,articlecolboismaerten} or in the setting of a sequence of manifold converging to some limit object \cite{articlecolboiscourtoisspectrum1form2, articleannecolbois, articlecolboiscourtois2, articlelottdifform2, articlelottspectrum,articlecolboismatei} or when degenerating as in \cite{artjammesvaleurpropres03}. Note also that some of its general properties were investigated in \cite{artlilowerbound,articledodziukdiffform,articlecolboiselsoufi,artjammesvaleurpropres07}, the first article of this list being also related to eigenvalue lower bounds (involving Sobolev's constants). Finally, let us mention the recent work \cite{articlelipstern1form} in which the authors give related results for the first eigenvalue of the Hodge Laplacian acting on 1-forms on hyperbolic 3-manifolds. \\

On a closed oriented Riemannian manifold, the associate Hodge Laplacian $\Delta_p =d \delta + \delta d$ acting on $p$-forms preserves the Hodge decomposition of the space of $L^2$-differential forms and its spectrum
decomposes accordingly in three types of eigenvalues. The eigenvalue $0$, corresponding to harmonic forms, whose eigenspace is identified to the $p$-th de Rham cohomology group $H^p_{\mathrm{DR}}(M, \RR)$, the exact spectrum, corresponding to exact $p$-eigenforms and the coexact spectrum, corresponding to coexact p-eigenforms.\\

In particular, as the spectrum of exact $1$-forms coincide with the non zero spectrum of functions, the study of the spectrum of $1$-forms reduces to the case of coexact $1$-forms and we will concentrate from now on 
on this case. We denote it by $$0 < \lambda^1_{1} \le \lambda^1_{2} \le ... \lambda^1_{k} \le ... .$$ 

 In order to motivate the introduction of our lower bound, one must start with a preliminary observation: the spectrum of the coexact 1-form is the same as the one of coexact (d-2)-forms. Indeed, the mapping that associates to a 1-form $\alpha$ the (d-2)-form $d^{-1} \star \alpha$, where $d^{-1}$ is the unique coexact primitive of $\star \alpha$, is by Hodge theory, an isomorphism. If one wants to study the 1-forms coexact spectrum, one can study equivalently the (d-2) coexact spectrum, which is indeed the perspective we have in mind in this article. \\

We define our 'Cheeger's constant' for the $1$-forms analogously to the isoperimetric constant. We observe that in the definition of the isoperimetric constant of $M$, one first removes 
a separating hypersurface $S$ in $M$ and get two connected components $V_1$ and $V_2$ of $M\backslash S$. This may be considered as increasing by $1$ the dimension of the $0$-th cohomology group
$H_{DR} ^0(M,\mathbb R)$, which precisely counts the number of connected components of $M$. Moreover, in the definition of the isoperimetric constant $$ h(M) := \inf_{ S \hookrightarrow M, \ [S] = 0} \ \frac{\mu_{d-1}(S)}{\min(\mu(V_1), \mu(V_2))},$$
the denominator may be interpreted as the minimum of the volumes of the two primitives of ${S}$. \\

In that spirit, we look at a way of increasing by $1$ the dimension of $H^{d-2}_{\DR}(M, \RR)$. We shall do it by cutting off a knot $\gamma$ in $M$, meaning by removing a smooth curve of $M$. Roughly speaking, the 
$(d-2)$-th absolute homology group of $M_{\gamma} := M \setminus \gamma$ is of dimension $\dim(H^{d-2}_{\DR}(M, \RR)) +1$ and is generated by the homology of $M$ together with 
a meridean $(d-2)$-sphere of the knot defined by $\gamma$. 
Following Cheeger's definition, we want to define the geometric cost of such a cutting by considering the ratio of the length of the curve $\gamma$ by the minimum 'area of a geometric primitive' of $\gamma$. 
%The most natural way to this area is perhaps as follows: 
For a real homologically trivial curve $\gamma$ of $M$ we set
\begin{equation}
	\label{eq.deminimalfcheeger1form}
	 A(\gamma) := \infi{ S_\gamma ,\,\,\partial S_{\gamma} = r \cdot \gamma} \frac{\area(S_{\gamma})}{r} \ ,
\end{equation}
where $r$ is any integer such that the curve $\gamma$ travelled up $r$ times, that we denoted $r \cdot \gamma$, is trivial in $H_1(M, \ZZ)$ and where the infinum ranges over all rectifiable current whose boundaries is $r \cdot \gamma$. Note that one can show, using Federer-Flemming compactness theorem for integral currents see \cite{artfedererflemming} or \cite[4.2.17]{livrefederergmt}, that $A(\gamma)$ is attained by a rectifiable current. \\

The integer $r$ is meant to encompass the case where $H_1(M, \ZZ)$ has a non trivial torsion part but will play no significant role. Note also that we could have chosen $r$ 
independently of $\gamma$, by setting, for example, $r$ as the less common multiple of all the possible orders of the torsion part of $H_1(M, \ZZ)$.

\begin{definition}
	\label{def.cheeger1forms}
	Let $M$ a closed Riemannian manifold. We set 
	$$ h^1 := \infi{\gamma \hookrightarrow M} \ \frac{l(\gamma)}{A(\gamma)} \ , $$
	where the infinum runs over all real homologically trivial smooth closed curves $\gamma$ and where $l(\gamma)$ is the length of $\gamma$.
\end{definition}

\begin{remark} \
\begin{itemize}
\item  The larger the class of 'surfaces' we consider when setting $A(\gamma)$, the better the constant $h^1$; it is meant to be a lower bound. We did not require for example that the chain is either connected or embedded. Ultimately, the class we want to work with is the largest one for which the area is well defined and for which one can use Stokes' theorem for smooth differential forms.
\item It is a priori not clear whether or not the above constant is positive. It seems well known to the expert but since we couldn't find a ready-to-use statement we devoted the appendix of this article to it. This proof was suggested to us by Laurent Mazet.
\end{itemize}
\end{remark}

Let us now state our main Theorem, which we think of as a Cheeger-like inequality for coexact 1-forms in the regime of small eigenvalues.
\begin{theoreme}[Main Theorem]
	\label{theo.cheegersurfaceminimale}
	\label{maintheo}
For any $a < b$ and any $D > 0$, there is a constant $C = C(a,b,D)$ such that for any closed orientable Riemannian manifold of diameter at most $D$ and sectional curvature bounded below by $a$ and above by $b$ we have 
	$$ \lambda^1_1 \ge \min(1, C \cdot h^1)^2 \ .$$
\end{theoreme}

{\bf Remarks.} 1. In the statement of the above theorem, in dimension 3, one can replace the two sided bounds assumption of the sectional curvature with the assumption that the Ricci curvature is bounded from below; the above inequality then holds in that case with a constant $C$ which only depends on a lower bound of the Ricci curvature and an upper bound of the diameter. We refer to Remark \ref{rem.dimension3} in Section \ref{sec.reductiontocurve} for more details. \\

2. In section \ref{section5}, we construct a family of metrics $(g'_\epsilon )_{0< \epsilon \leq 1}$ on the $3$-sphere with uniformly bounded by below sectional curvature, diameter tending to $\infty$ when $\epsilon$ tends to $0$, such that $\lim _{\epsilon \to 0} \lambda ^1 _1 (S^3, g'_\epsilon) =0$ and $h^1 (S^3, g'_\epsilon) \geq C'>0$. This example shows that in
Theorem \ref{theo.cheegersurfaceminimale},
the dependency of the constant $C$ in the diameter is necessary. Since the family of metric $g_{\epsilon}$ also carries negative curvature, one can scale the metrics to have diameter 1 but Ricci curvature unbounded from below, which also shows the assumption on the Ricci curvature are necessary. \\

3. Our definition of $h^1$ is very close in spirit to the definition of Cheeger's constant $h^0$. However, the statement of our main theorem and the counterexample of Section \ref{sec.counterexample} show a different feature for the 1-forms case: the lower bound is not universal since it depends on the global geometry, through the diameter, and the local one, through the Ricci curvature. One can wonder if the definition of 'least area of a geometric primitive' may be relaxed to get a lower bound as in Theorem \ref{theo.cheegersurfaceminimale} but with an universal constant. Such a question will be addressed in a forthcoming article, see Section \ref{sec.remarks}. \\

4. There is a noticeable difference between small eigenvalues for the functions and for the differential forms. The geometry of $M$ gives more constraints in the first case.  Indeed, on the one hand, as shown in \cite{Gro} and \cite{BBG}, a lower bound on the Ricci curvature together with an upper bound on the diameter of $M$ prevent the Cheeger's constant -and thus, by Cheeger's Theorem, $\lambda _1 (M)$- to be small.
On the other hand, the same assumptions on Ricci curvature and diameter do not prevent from small eigenvalues for differential forms. For example, the family of collapsing Berger's metrics on the $3$-sphere $(S^3, g_\epsilon)_{0< \epsilon \leq 1}$
have bounded sectional curvature and diameter while the first coexact eigenvalue for the $1$-forms satisfies $\lambda ^1 _1 (M, g_\epsilon) = 4 \epsilon ^2$, cf.
\cite{articlecolboiscourtoisspectrum1form}. \\

5. In section \ref{sec.example}, we show that the inequality of Theorem \ref{theo.cheegersurfaceminimale} is sharp up to a constant: indeed, the Berger's spheres satisfy $\lambda ^1 _1(S^3,g_\epsilon) = 4 \epsilon ^2$ and $h^1 (S^3, g_\epsilon) \sim a\epsilon$ for some constant $a>0$ when $\epsilon$ tends to $0$. \\

6. In dimension 3 and 4, the full spectrum of $M$ (the union of the spectrum of all degrees) corresponds to the 1-forms spectrum (exact and coexact). Theorem \ref{maintheo} together with Cheeger's inequality give a lower bound on the full spectrum of $M$ in these dimensions. \\

\textbf{A word on the proof.} Since the two next sections of this article readily address the proof of Theorem \ref{maintheo}, we will not discuss it here in detail. We however shall quickly explain how the constant $h^1$ and $\lambda^1_1$ may sensibly be related. The key for proving the classical Cheeger's inequality is to relate the geometry of the level sets of a function with its Rayleigh quotient. In the function case, a somewhat natural tool is the coarea formula, which is indeed the starting point of Cheeger's proof. The main difficulty, as noticed by P. Buser, to generalise Cheeger's inequality to the setting of differential forms is to replace the use of the coarea formula by another tool. Let us explain what is this tool in our context, meaning in the setting of coexact (d-2)-forms. \\ 

Let then $\alpha$ be such a differential forms. We set $X := (\star d \alpha)^{\sharp}$, which is a vector field. The fact that $d \alpha$ is exact implies that the flow of $X$ preserves the Riemannian measure $\mu$. This measure preserving property implies that the following identity holds for any 1-form $\beta$ and any $T>0$:
 \begin{equation}
	\label{eq.coarea1form}
	 \inte{M} d \alpha \wedge \beta  = \inte{M} \frac{1}{T} \inte{\gamma(x,T)} \beta \ d \mu(x) \ , 
\end{equation}
where $\gamma(x,T)$ is the piece of trajectory of the $X$-flow issued from $x$ up to time $T$. This formula may be read by duality: the (d-1)-form $d \alpha$ is dual to the foliation induced by the flow of $X$, or in other words, dual to the trajectories of $X$. Instead of directly relating $h^1$ to $\lambda_1^1$, we shall first relate $h^1$ to the $L^2$-Rayleigh quotient of $\alpha$ (as in Cheeger's proof). The starting point of the analysis is to use the above formula with the 1-form $\beta$ such as its integral over $\gamma(x,T)$ is the length of $\gamma(x,T)$. The left member will then turn out to be the $L^1$ norm of $d \alpha$, the numerator of the $L^1$-Rayleigh quotient. Note that this step is completely analogous to Cheeger's starting point. In order to insert the constant $h^1$, one must work with homologically trivial closed curves, which may not be the case of the $\gamma(x,T)$'s. Therefore, instead of using the $\mu$-average over all trajectories as in the right hand side of Formula \ref{eq.coarea1form}, we will produce 
a family of 'closed trajectories' $\Gamma(n,T) _{n\in \mathbb N, \, T>0}$ which are homologically trivial and such that the integration currents $$ \frac{1}{nT} \inte{\Gamma(n,T)} \cdot$$ acting on $1$-forms approximate, when $T, \,n \to \infty$, the term
$$
 \inte{M} d \alpha \wedge \cdot
$$
of the left hand side of Formula  \ref{eq.coarea1form}. These curves $\Gamma(n,T) _{n\in \mathbb N, \, T>0}$ will be produced by a probalistic argument,
choosing at random $n$ trajectories of the flow of $X$ up to time $T$ that are glued altogether in a careful way in order to get te two above properties. This point, the heart of this article, will be discussed in details in Section \ref{sec.proofcurve}. \\

 Let us conclude by reviewing the historical works that inspired us. The starting point is Arnold's idea to use the geometry of the trajectories of a flow (the self linking number or isoperimetrical ratio) in order to interpret an analytical quantity associated to the dual 1-form (the helicity or its energy). Arnold's idea may actually be rooted in Schwartzmann's theory of asymptotic cycles \cite{articleschwartzman}, generalised by Sullivan \cite{articlesullivancycles}, which relates the average homology class defined by the orbits of a measure preserving flow to the so called `asymptotic cycle' associated to the pair (measure/vector field). It is on this rooting idea that relies the key step of the proof: closing up a random trajectories for a small cost in length in such a way that the resulting curve becomes homologically trivial. \\
 
\textbf{Acknowledgement.} We would like to thank Laurent Mazet for having indicated to us the proof of our appendix statement. We would also like to thank Antoine Julia for his explanations of Federer-Fleming compactness theorem as well as Julien Marché for his explanations of knots complementary topology. We finally would like to thank Sylvestre Gallot for having communicated to us the paper \cite{Gal2}.

\section{Proof of our main theorem; reduction to Proposition \ref{propexistencecurve}}
\label{sec.reductiontocurve}

This section is devoted to the first step of the proof of our main theorem. The next subsection aims at recalling some basic facts together with introducing the notations we will use.

\subsection{Generalities and notations}
\label{subsec.generalities}
Let $M$ be an oriented Riemanian manifold of dimension $d$. For any $x \in M$ we denote by $T_x M$ of tangent plane at $x$, by $\left< \cdot, \cdot \right>_x$ the metric on $ T_x M$, by $TM$ the tangent bundle of $M$ and by $\left< \cdot, \cdot \right>$ the Riemannian metric on $M$. If $X_x \in T_x M$ we denote by $|X_x|$ is norm with respect to the metric of $T_x M$. If $X$ is a vector field, we denote by $|X|$ the function $x \mapsto |X_x|$. If $M$ has dimension $d$ and $0\leq p \leq d$, we denote by $\Lambda_x^p(M)$ the set of $p$ alternating multilinear form on $T_x(M)$. A differential $p$-form on $M$ is a section of the fibre bundle $\Lambda^p(M)$. \\

If $\alpha \in \Lambda^1(M)$ we denote by $\alpha^{\sharp}$ its musical dual vector field $X$ defined as
	$$ \alpha(Y) = \left<X, Y \right> \ ,$$
for any vector field $Y$. Conversely, if $X \in T_x(M)$, we denote by $X^{\flat}$ the 1-form $\alpha$ defined as above. This duality is usually refered as the musical duality. The above duality endows the vectorial spaces $\Lambda_x^1(M)$ with a scalar product and a norm as well. We shall keep denoting them by $\left< \cdot, \cdot \right>$ and $| \cdot |$. There is another way to think of this scalar product thanks to the Hodge star operator, or Hodge operator for short, which is defined as follows. 
\begin{definition}
Let $M$ be an oriented $d$-dimensional Riemannian manifold, $0 \le p \le n$ and $\alpha \in \Lambda^p(M)$. The Hodge operator, denoted by $\star$, is defined as the unique operator from $\Lambda^p(M)$ to $\Lambda^{n-p}(M)$ which satisfies the following: for any $x \in M$ and for any oriented orthornormal basis $(\alpha_1, ..., \alpha_n)$ of $\Lambda_x^1(M)$ 
		$$ \star (\alpha_1 \wedge ... \wedge \alpha_p) = \alpha_{p+1} \wedge ... \wedge \alpha_n \ . $$
\end{definition}

One can easily check from the above definition that $\star$ is an involution, up to sign:
	\begin{equation}
		\label{eqhodgeinvolutive}
			 \star \star = (-1)^{p(n-p)} \Id \ .
	\end{equation}
Note also that the above definition readily implies 
	$$  \star \ \mu = 1 \ ,$$
where we denoted by $ \mu$ the Riemannian volume form associated to the metric. We shall also refer to the Riemannian measure as $\mu$, slightly abusing the notation. The scalar product induced by the musical isomorphism is easily generalisable by using the Hodge operator. We define the punctual scalar product of two $p$-forms $\alpha, \beta$ as the following real valued function
	\begin{equation}
		\label{eqhodgescalarproduct}
		 (\alpha, \beta) := \star (\alpha \wedge \star \beta) \ .
	\end{equation}
	Note that Equation \eqref{eqhodgeinvolutive} is mandatory to reach the symmetry required for a scalar product.  This symmetry implies in particular 
	$$ | \star \beta | = | \beta | \ . $$
The norm induced by the previous scalar product extends the one already mentioned on 1-forms; which is why we will also denote it by $| \cdot |$. The punctual norm can be integrated to endow $\Lambda^p(M)$ with $L^k$ norms for $ 1 \le k \le \infty$ . For any $\beta \in \Lambda^p(M)$ we define 
 $$ || \beta ||_{k} := \left( \inte{M} \ |\beta|^k \ d  \mu \right)^{1/k} \ . $$
 which makes (for any $k \ge 0$ and any $0 \le p \le n$) the Hodge operator an isometry from $\Lambda^p(M)$ to $\Lambda^{n-p}(M)$, both of them endowed with the $L^k$ norm. In the case where $k =2$, the space of $L^2$-differential forms comes with a global scalar product. Indeed, for $\alpha, \beta \in \Lambda^p(M)$ we set
		$$ \left< \alpha, \beta \right> := \inte{M} \alpha \wedge \star \beta \ , $$
slightly abusing the notation since we already use the brackets for the Riemannian metric. The exterior derivative $d$ has an adjoint $\delta$ with respect to the above scalar product: 
	$$ \delta := (-1)^{n(p-1) +1} \star d \star \ . $$
Note that $\delta$ maps $\Lambda^p(M)$ to $\Lambda^{p-1}(M)$. Stokes formula then shows that for any $\alpha \in \Lambda^{p-1}(M)$ and any $\beta \in \Lambda^{n-p}(M)$ we have
	$$ \left< d \alpha, \beta \right> = \left< \alpha, \delta \beta \right> \ .$$
We now have the material to define the Hodge Laplacian.

\begin{definition}
Let $M$ be a closed oriented Riemannian manifold of dimension $d$. For $0 \le p \le d$ we define de Hodge Laplacian $\Delta_p$ which acts on $\Lambda^p(M)$ as 
	$$ \Delta_p := d \delta + \delta d \ . $$
\end{definition}

By construction, this operator is symmetric, meaning that for any $\alpha, \beta \in \Lambda^p(M)$ we have 
	$$ \left< \Delta_p(\alpha) , \beta \right> = \left< \alpha ,\Delta_p(\beta) \right> \ . $$
	
	The two following simple lemma will be used in the sequel.
	
\begin{lemma}
	\label{lemme.interiorproduct}
	Let $M$ be a closed oriented Riemannian manifold of dimension $d$. For any 1-form $\beta$ and any vector field $Y$ we have  
	\begin{equation}
			  \ i_Y(\beta) = \star (Y^{\flat} \wedge \star \beta ) =   (Y^{\flat},\beta)\ .
	\end{equation}
\end{lemma}

The proof of the above identity consists in evaluating the above formula on an oriented orthonormal basis of the form
	$$ \left( \frac{Y}{||Y||}, X_2, ... , X_n \right) $$ 	
at each point of $M$ where $Y$ does not vanish and is left to the reader. We shall also need the notion of mass of a differential form.
	
\begin{definition}	
	\label{def.mass}
	Let $M$ be a Riemannian manifold of dimension $d$. for any $p$-form $\beta$ we define its punctual mass as
	\begin{equation}
		\label{eqpunctualnormpforms}
			 \pmass(\beta)(x) := \supr{X_1, ..., X_p} \ \beta(X_1, ..., X_p) \ ,
	\end{equation}
where the sup runs over the set of orthonormal systems $(X_1,..., X_p)$ of vectors in $T_x M$. We also define its mass as 
	$$ \mass(\beta) := || \pmass(\beta) ||_{\infty} \ . $$
\end{definition}

The punctual norm of a differential form and its mass are closedy related as testified by the conclusion of the following lemma whose proof is also left to the reader. 

\begin{lemma}
	\label{lemma.massinfty}
	Let $M$ be a Riemannian manifold of dimension $d$. Then, for any $p$-form $\beta$ of $M$ we have 
		$$ \mass(\beta) \le || \beta ||_{\infty} \ .$$ 
\end{lemma}
  	
\subsection{Reduction to Proposition \ref{propexistencecurve}}
\label{subsec.reductiontocurve}

This subsection aims at reducing the proof of our main theorem to the key proposition stated below, whose proof will occupy entirely the section \ref{sec.proofcurve}.
Given a coexact form $\alpha \in \Lambda ^{d-2} (M)$ and $\beta$ any 1-form such that $\star \alpha = d \beta$,
the Proposition provides a family of homologically trivial closed curves $\Gamma (n, T) _{n \in \mathbb N , T>0}$ which, roughly speaking, can be thought of as 
typical long trajectories (up to time $T$) of the vector field $X$ dual to $\star d \alpha$ glued up altogether.  The main properties of these curves are that $ || d \alpha||_{1} $ is approximated by 
	$$\frac{ l(\Gamma(n,T))}{nT}$$ 
	and 
$ || \alpha||_2^2$ is approximated by
 $$\frac{1}{nT} \inte{\Gamma(n,T)} \beta \ ,$$ 
 leading to the relation between the Rayleigh quotient of $\alpha$ with the constant $h^1$. We now state this Proposition and prove how it implies Theorem \ref{maintheo}.

\begin{proposition}
	\label{propexistencecurve}
	Let $\alpha \in \Lambda^{d-2}(M)$ be coexact and $\beta$ any 1-form such that $\star \alpha = d \beta$. For all $\epsilon > 0$ there is a time $T > 0$, an integer $n \in \NN$ and a closed curve $\Gamma(n,T)$ such that the following holds simultaneously,
\begin{enumerate}
\item the closed curve $\Gamma(n,T)$ writes as an union 
	$$ \Gamma(n,T) = \gamma(n,T) \cup \nu(n,T) \ , $$
with $$ \frac{l(\nu(n,T))}{nT}  \le \epsilon \ ;$$
\item the curve $\gamma(n,T)$ satisfies that
$$	\Big| || d \alpha||_{1} - \frac{ l(\gamma(n,T))}{nT} \Big| \le \epsilon \ ;  $$
\item the curve $\Gamma(n,T)$ is homologically trivial as a cycle with real coefficients;
\item the curve $\gamma(n,T)$ satisfies that
$$	\Big| || \alpha||_2^2 - \frac{1}{nT} \inte{\gamma(n,T)} \beta \Big| \le \epsilon \ .  $$
\end{enumerate}
\end{proposition}

As it has been already alluded to, the role of the curves $\Gamma(n,T)$ is to get the analogue of Formula \ref{eq.coarea1form} saying that
 the integration currents on $1$-forms 
 $$ \frac{1}{nT} \inte{\Gamma(n,T)} \cdot \ \ \text{ approximate } \ \ 
 \inte{M} d \alpha \wedge \cdot \ \ \text{ when } \ \ T, n \to \infty \ . $$ Applying this successively to appropriate $1$-forms give the items (2), (3), and (4), see section \ref{sec.proofcurve}.

\bigskip
\textbf{Proof of (Proposition \ref{propexistencecurve} $\Rightarrow$ Theorem \ref{maintheo}).} In order to clarify the structure, we split the proof in two main steps, the first one being contained in 

\begin{lemma}
\label{lem.keyproptomaintheo1}
Let $\alpha \in \Lambda^{d-2}(M)$ be coexact. We have
	$$  \frac{||d \alpha||_{2}}{|| \alpha ||_2} \cdot \left(  \frac{\mu(M)^{\frac{1}{2}} \cdot ||  \alpha||_{\infty} }{||\alpha||_2} \right) \ge  h^1 \ . $$
\end{lemma}

\textbf{Proof.}  Let $\alpha$ as above and $\beta$ any 1-form such that $\star \alpha = d \beta$. The proof of the above lemma uses all the items of Proposition \ref{propexistencecurve} one by one. The first item of \ref{propexistencecurve} readily gives that for all $\epsilon > 0$ there exist $n \in \mathbb N, T>0$, curves $\Gamma(n,T)$ and $\gamma(n,T)$ such that 
$$ \frac{l(\gamma(n,T))}{nT} \ge \frac{l(\Gamma(n,T))}{nT}  - \epsilon \ .$$
Combined with the second item of Proposition \ref{propexistencecurve} we then get
	\begin{align*}
		||d \alpha ||_{1} & \ge \frac{l(\gamma(n,T))}{nT}  - \epsilon  \\
			& \ge \frac{l(\Gamma(n,T))}{nT}  - 2 \epsilon  \ .
	\end{align*} 

Since $\Gamma(n,T)$ is homologically trivial as a real cycle, there is an integer $r$ such that $ r \cdot \Gamma(n,T)$ bounds a 2-chain $S_{\Gamma(n,T)}$ and which realises $A(\gamma)$. Our definition of $h^1$ then implies that 
	\begin{align*}
		 l(\Gamma(n,T))  & \ge h^1 \cdot \frac{\area(S_{\Gamma(n,T)})}{r} \\
		 	 & \ge \frac{h^1}{r} \  \supr{ \mass(\omega) \le 1}  \ \inte{S_{\Gamma(n,T)}} \ \omega \\
		 	 & \ge \frac{h^1}{r} \ \supr{ || \omega ||_{\infty} \le 1}  \ \inte{S_{\Gamma(n,T)}} \ \omega \ ,
	\end{align*}
because of Lemma \ref{lemma.massinfty}. In particular, one has
$$ \supr{ || \omega ||_{\infty} \le 1}  \ \inte{S_{\Gamma(n,T)}} \ \omega \ \ge \inte{S_{\Gamma(n,T)}} \ \frac{\star \alpha}{|| \star \alpha||_{\infty} } = \frac{1}{||\alpha||_{\infty} } \  \inte{S_{\Gamma(n,T)}} \ d \beta \ , $$
since $\star$ is an isometry and by definition of $\beta$. Stokes theorem then yields 
$$	\supr{ || \omega ||_{\infty} \le 1}  \ \inte{S_{\Gamma(n,T)}} \ \omega \ \ge \frac{1}{||\alpha||_{\infty} } \  \inte{ r \cdot \Gamma(n,T)} \ \beta = \frac{r}{||\alpha||_{\infty} } \  \inte{ \Gamma(n,T)} \ \beta \ , $$
and then 
	\begin{equation}
	\label{eq.demmaintheo100}
		||d \alpha ||_{1}  \ge \left(  \frac{h^1}{||  \alpha||_{\infty} } \ \frac{1}{nT} \ \inte{\Gamma(n,T)} \ \beta \right) -  2 \epsilon \ .		
	\end{equation}
We chop off the above integral in two pieces with the objective of using point (4) of Proposition \ref{propexistencecurve}. 
$$ \frac{1}{nT} \ \inte{\Gamma(n,T)} \ \beta = \frac{1}{nT} \ \inte{\gamma(n,T)} \ \beta + \frac{1}{nT} \ \inte{\nu(n,T)} \ \beta \ . $$
Because of Proposition \ref{propexistencecurve} (1) we have 
$$ \Big| \frac{1}{nT}  \inte{\nu(n,T)} \ \beta \Big| \le || \beta||_{\infty} \cdot \epsilon \ ,$$
and then
$$  \ \frac{1}{nT} \ \inte{\Gamma(n,T)} \ \beta \ge \ \frac{1}{nT} \ \inte{\gamma(n,T)} \ \beta -  || \beta ||_{\infty} \cdot \epsilon \ .$$
Combined with \eqref{eq.demmaintheo100}, it yields 

	\begin{align*}
		||d \alpha ||_{1} & \ge \left( \frac{h^1}{||  \alpha||_{\infty} } \ \frac{1}{nT} \ \inte{\gamma(n,T)} \ \beta \right) - \frac{ h^1 \cdot || \beta ||_{\infty}}{|| \alpha||_{\infty}} \epsilon - 2 \epsilon \ .		
	\end{align*} 
	
We now use Proposition \ref{propexistencecurve} (4) to get 
\begin{align*}
	 \ \frac{1}{nT} \ \inte{\gamma(n,T)} \ \beta \ & \ge || \alpha||^2_{2} - \epsilon  \ .
\end{align*}
Therefore,

\begin{align*}
		||d \alpha ||_{1} & \ge \left( \frac{h^1}{||  \alpha||_{\infty} } \  ( || \alpha||^2_{2} - \epsilon ) \right) - \frac{ h^1 \cdot || \beta ||_{\infty}}{|| \alpha||_{\infty}} \epsilon - 2 \epsilon \\ 
		& \ge  \frac{h^1 \cdot || \alpha||^2_{2}  }{||  \alpha||_{\infty} } \  - C \epsilon \ , 
	\end{align*} 	
with $C= C(\alpha, \beta, h^1)$ some constant independent of $\epsilon$. Since none of the above quantities depends on $\epsilon$ (except $\epsilon$ itself, of course) one can let $\epsilon \to 0$ and get
	$$  	||d \alpha ||_{1} 	\ge  \frac{h^1 \cdot || \alpha||^2_{2}  }{||  \alpha||_{\infty} } \ ,  $$
We use the Cauchy Schwartz inequality to give an upper bound to the above left $L^1$-norm and get
	$$  	||d \alpha||_{2} \ \mu(M)^{\frac{1}{2}} \ge  \frac{h^1 \cdot || \alpha||^2_{2}  }{||  \alpha||_{\infty} } \  $$
which yields to 
	$$  \frac{||d \alpha||_{2}}{|| \alpha ||_2} \cdot \left(  \frac{\mu(M)^{\frac{1}{2}} \cdot ||  \alpha||_{\infty} }{||\alpha||_2} \right) \ge  h^1 \ , $$
which is the desired inequality. \hfill $\blacksquare$ \\

Let then $\beta_0 \in \Lambda^1(M)$ be the eigenform of $\Delta_1$ associated to $\lambda^1_1$. As already emphasized, the Hodge theory guarantees that there is an unique coexact (d-2)-forms $\alpha_0 \in \Lambda^{d-2}(M)$ defined as $d \alpha_0 = \star \beta_0$ which has same eigenvalue $\lambda^1_1$. In particular, we have 
$$ \frac{||d \alpha_0||_{2}}{|| \alpha_0 ||_2} = \sqrt{\lambda^1_1} \ ,$$

and then, by Lemma \ref{lem.keyproptomaintheo1},
\begin{equation}\label{funda-ineq}
  \sqrt{\lambda^1_1} \cdot \left(  \frac{\mu(M)^{\frac{1}{2}} \cdot ||  \alpha_0||_{\infty} }{||\alpha_0||_2} \right) \ge  h^1 \ . 
 \end{equation}

Our Main Theorem \ref{theo.cheegersurfaceminimale} now follows from the above (\ref{funda-ineq}) and a classical upper bounds of the ratio  $$\frac{\mu(M)^{\frac{1}{2}} \cdot ||  \alpha_0||_{\infty} }{||\alpha_0||_2} \ ,$$
for any eigenform $\alpha _0$ of the Hodge Laplacian. This upper bound depends on the associated eigenvalue and on an upper bound on the diameter and a lower bound on the curvature operator of the underlying manifold $M$. Since we could not find any ready-to-quote reference, we decided to encompass this classical result within the
\begin{proposition}
\label{prop.keyproptomaintheo1}
For any $a < b$ and, $D > 0$ and any $ 0 \le \lambda \le 1$, there is a constant $C_1 = C_1(a,b,D)$ such that for any closed orientable Riemannian manifold of diameter at most $D$ and sectional curvature bounded below by $a$ and above by $b$ and for any p-form $\alpha$ such that $\Delta \alpha = \lambda \alpha$ then 
$$ \frac{\mu(M)^{\frac{1}{2}} \cdot ||  \alpha_0||_{\infty} }{||\alpha_0||_2} \le C_1 \ . $$
\end{proposition}

\textbf{Proof.} The key tool of the proof is the following
 \begin{theoreme} \cite{Gal1,Gal2},{\cite[page 391]{artberardbochner}}
\label{theo.comparisonnormeigensections}
For any $R_{\min} \in \RR$, any $D > 0$ andany $B>0$, there is a constant $C_1 = C_1(R_{\min},D, B)$ such that for any closed orientable Riemannian manifold of diameter at most $D$ and Ricci curvature bounded below by $R_{min}$ and for any p-form $\alpha$ which satisfies
	$$ \Delta |\alpha| \le B |\alpha|$$ 
in the distributional sense, then
	$$ \frac{ \mu(M)^{1/2} \cdot || \alpha||_{\infty} }{ || \alpha||_2 } \le C_1 \ .$$
\end{theoreme}
The above theorem is actually stronger since it can be stated for sections of a general fibre bundle. Moreover, the constant $C_1$ may be exhibited as the first zero of an explicit differential equation with parameters explicitly related to the constants $D$, $R_{\min}$ and $B$. The $\min$ appearing in the lower bound of our Theorem \ref{maintheo} can now be explained: if $\lambda^1_1 \ge 1$ there is nothing to prove so that one can assume that $\lambda^1_1 \le 1$. Under this extra assumption, in order to use the above Theorem \ref{theo.comparisonnormeigensections}, one is left to show the

\begin{lemma}
\label{lemma.norm.of.eigenform}
For any $a < b$ and any $D > 0$, there is a constant $B = B(a,b,D)$ such that for any closed orientable Riemannian manifold of diameter at most $D$ and sectional curvature bounded below by $a$ and above by $b$ and for any p-form $\alpha$ such that $\Delta \alpha = \lambda \alpha$ then 
	$$ \Delta |\alpha| \le B |\alpha|$$ 
holds in the distributional sense. 
\end{lemma}
\textbf{Proof.} The proof is classical and can be extracted from \cite{Gal2} or \cite{artberardbochner} as well. We write it down for the reader's convenience. It relies on the second Kato's inequality which asserts that the following inequality holds in the distributional sense
\begin{equation}
	\label{eq.kato}
	 |\alpha| \Delta | \alpha| \le \left< \bar{\Delta}_p \alpha, \alpha \right> \ , 
\end{equation}
where $\bar{\Delta}_p$ is the Bochner Laplacian acting on p-forms, also called rough Laplacian. This Laplacian is related to the Hodge one by the Weitzenböck formula,
	$$ \bar{\Delta}_p \alpha = \Delta_p  \alpha - \mathrm{Ric}_p(\alpha) \ , $$
where $\mathrm{Ric}_p$ stands for the curvature operator acting on p-forms. This operator is also called Ricci curvature operator, justifying the notation, since it coincides with the Ricci operator in restriction on 1-forms. We shall use the later terminology to avoid confusion with the classical curvature operator. Therefore,
\begin{equation}
\label{eq.lowerboundriccurvop}
	 \left< \bar{\Delta}_p \alpha, \alpha \right> = \left< \Delta_p  \alpha, \alpha \right> -  \left< \mathrm{Ric}_p(\alpha), \alpha \right> \ .
\end{equation}

In order to use Theorem \ref{theo.comparisonnormeigensections}, one must then bound from below the term $\left< \mathrm{Ric}_p(\alpha), \alpha \right>$. On the one hand, one can retrieve an eigenvalue lower bound for the Ricci curvature operator acting on p-forms from a lower bound $\mathrm{\lambda}_{\min}$ on the eigenvalues of the curvature operator \cite[page 264]{artgallotmeyeropcourbure};
	$$ \left< \mathrm{Ric}_p(\alpha), \alpha \right> \ge p(n-p) \mathrm{\lambda}_{\min} |\alpha|^2 \ . $$
	
On the other hand, it is well known, see for example \cite[Proposition 3.8]{artbourgignonopcourbure}, that the curvature operator has eigenvalues bounded from below (actually from above as well) provided that the sectional curvature is bounded above and below. In particular there is a constant $C = C(a,b)$ (recall that $a$ and $b$ are the bounds on the sectional curvature) such that 	
	$$  \left< \mathrm{Ric}_p(\alpha), \alpha \right> \ge C \cdot | \alpha |^2 \ . $$

Since we assume $\Delta \alpha = \lambda \alpha$ with $\lambda \le 1$ we have
\begin{align*}
\left< \bar{\Delta}_p \alpha, \alpha \right> & \le |\alpha|^2 + \max(- C,0) \cdot  |\alpha|^2  \\
	& \le (1 + \max(- C,0)) \cdot |\alpha|^2 
\end{align*}

Together with \eqref{eq.kato} it yields 
		$$  |\alpha| \Delta | \alpha| \le (1 + \max(- C,0)) \cdot |\alpha|^2   $$ 
in the distributional sense. The above inequality holds at points where $\alpha = 0$: these points are local minimum of $|\alpha|$ which implies that $\Delta |\alpha| \le 0$ (in the distributional sense) at such points. Dividing both members by $|\alpha|$ whenever $\alpha \neq 0$ yields to the desired conclusion by setting $B := (1 + \max(- C,0))$.  \hfill $\blacksquare$ $\blacksquare$ \\

The Main Theorem \ref{maintheo} is therefore a consequence of Proposition \ref{prop.keyproptomaintheo1} and Lemma \ref{lem.keyproptomaintheo1} by setting
$C:= C_1 ^{-1}$.
 \hfill $\blacksquare$ \\

\begin{remark} 
\label{rem.dimension3} 1. In the above demonstration we only use a lower bound of the Ricci curvature (in order to apply Theorem \ref{theo.comparisonnormeigensections}) and  
a lower bound of the Ricci curvature operator acting on (d-2)-forms (in the proof of Lemma \ref{lemma.norm.of.eigenform}). In dimension 3, one has $d-2 =1$ and the Ricci curvature operator acting on 1-form coincides with the Ricci curvature, therefore a lower bound on the Ricci curvature is enough to prove our theorem. \\

2. Let us clarify the assumptions we need in our main theorem in the light of the proof: the two sided bounds of the sectional curvature can actually be replace with a lower bound on the eigenvalues of the curvature operator. Indeed, Inequality \eqref{eq.lowerboundriccurvop} shows that such a lower bound implies a lower bound on all the Ricci curvature operators acting on forms of any degree, in particular a lower bound on the Ricci curvature (see item 1 of this remark). In particular again, one can also use directly Inequality \ref{eq.lowerboundriccurvop} to show that a lower bound on the eigenvalues of the curvature operator gives a lower bound on the Ricci curvature operator acting on (d-2)-forms, allowing one to conclude the proof of Lemma \ref{lem.keyproptomaintheo1}. 
\end{remark}

\section{Proof of Proposition \ref{propexistencecurve}}
\label{sec.proofcurve}

It remains then to show that Proposition \ref{propexistencecurve} holds, proposition that we recall here with a slightly different but equivalent statement: the only difference resides in the new constant $C > 0$ which we introduce in order not to discuss so much the $\epsilon$ along the proof.

\begin{proposition*}\label{2-10*}
	Let $\alpha \in \Lambda^{d-2}(M)$ be coexact and $\beta$ any 1-form such that $\star \alpha = d \beta$. There is a constant $C > 0$ such that for all $\epsilon > 0$ there is a time $T > 0$, an integer $n \in \NN$ and a closed curve $\Gamma(n,T)$ such that the following holds simultaneously,
\begin{enumerate}
\item the closed curve $\Gamma(n,T)$ writes as an union 
	$$ \Gamma(n,T) = \gamma(n,T) \cup \nu(n,T) \ , $$
with $$ \frac{l(\nu(n,T))}{nT}  \le C \cdot \epsilon \ ;$$
\item the curve $\gamma(n,T)$ satisfies that
$$	\Big| || d \alpha||_{1} - \frac{ l(\gamma(n,T))}{nT} \Big| \le \epsilon \ ;  $$
\item the curve $\Gamma(n,T)$ is homologically trivial ;
\item the curve $\gamma(n,T)$ satisfies that
$$	\Big| || \alpha||_2^2 - \frac{1}{nT} \inte{\gamma(n,T)} \beta \Big| \le \epsilon \ .  $$
\end{enumerate}
\end{proposition*}

Beforing entering into the core of the proof, we sketch it and comment, again, the statement. Recall that for $X$ the vector field dual to $\star d \alpha$ and for $x\in M$, the trajectory of the $X$-flow starting at $x$ up to time $T$ have been denoted
by $\gamma (x,T)$.
Let us describe first how the curve $\gamma(n,T)$ is contructed: it consists of the union of $n$ pieces of trajectory $(\gamma(x_i, T))_{1 \le i \le n}$ starting at some points $x_1, ... , x_n \in M$. Note that, given $T > 0$, $\gamma(n,T)$ is completely determined by the starting points $x_1, ..., x_n$. The points $x_1, ...x_n$ won't be explicitely constructed but we will show their existence by a probabilistic argument: we will take $n$ random points independently with respect to the Riemannian measure $\mu$. The curve $\gamma(n,T)$ has now to be thought of as a random curve. The key remark is that the random integration current (defined on 1-forms)
	$$ \frac{1}{nT} \inte{\gamma(n,T)} \cdot $$
satisfies the law of large number and must behave, when $n$ tends to infinity, as its mean 
$$
\frac{1}{T} \int _M \left(\inte{\gamma (x,T)} \cdot \right) d\mu (x)
$$
which in turn can easily be identified to 
 		$$ \inte{M} d \alpha \wedge \cdot \ , $$ 
by Formula \ref{eq.coarea1form}. From this remark, the Points (2) and (4) of Proposition \ref{2-10*} will then boils down in evaluating appropriately chosen forms to the integral currents. \\

Let us now discuss the points (1) and (3) which are non trivially related and depend on the curve $\nu(n,T)$. We refer to the curve $\nu(n,T)$ as the closing-with curve: its only purpose is to close up the curve $\gamma(n,T)$ with a 'relatively small cost' in length (as in point (1)) but in a way that $\Gamma(n,T)$ is closed and homologically trivial (as in point (3)). Following Arnold's steps, we first glue altogether the $n$ trajectories of $\gamma(n,T)$ with uniformly bounded segments. The resulting curve $\gamma_c(n,T)$ is then close but has no reason to be homologically trivial. In order to make it so, we shall 'unroll' it homologically with the help of a system of curves generating the homology of $M$. The cost in length of such an unrolling must be controlled, meaning that one must know \textit{a priori} that the homology class of $\gamma_c(n,T)$ is relatively small with respect to its length. This points will also be handled by using the law of large number, thinking of $\gamma(n,T)$ as a random cycle, in the spirit of Schwartzman's theory of asymptotic cycles. \\

\begin{figure}[h]
\begin{center}
	\def\svgwidth{0.8 \columnwidth}
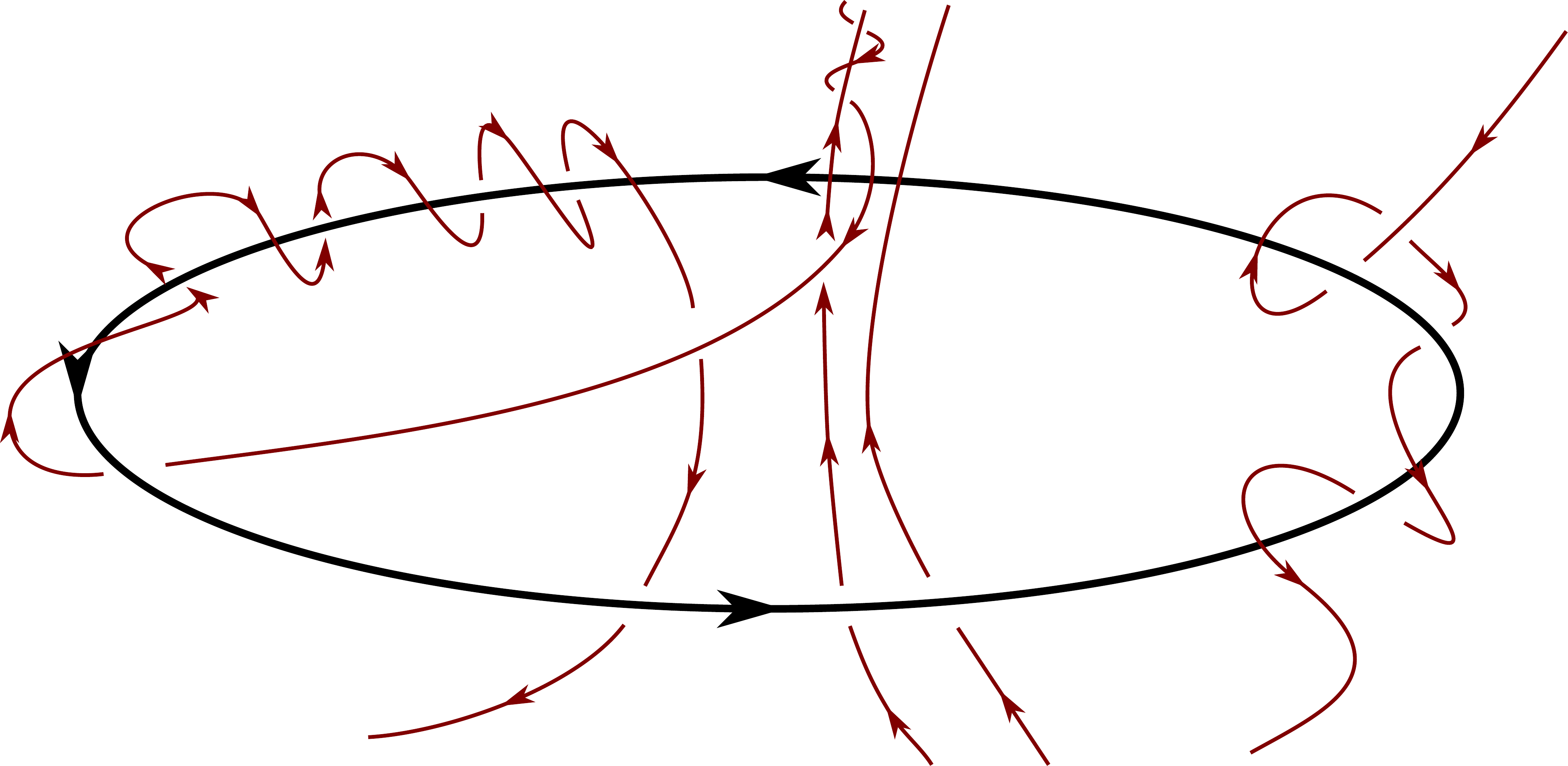
	\end{center}
\caption{This picture shows a typical realisation of a curve $\gamma_c(n,T)$ (in red) in the manifold $M = \RR^3 \setminus \mathbb{S}^1$ whose homology is generated by a meridian of the black circle. The red curve tends to be dense in $M = \RR^3 \setminus \mathbb{S}^1$ and is not necessary homologically trivial. However, it tends to link with the black circle (meaning not to be homologically trivial) not too often (once here) compared to its length.}
	\label{fig.curvestat}
\end{figure}

In order to emphasize the dependences between all the quantity previously introduced (which is important because of how is $T$ set) we shall first detail the construction of the closing-with curve $\nu(n,T)$. 

\subsection{Step 1: the closing-with curve $\nu(n,T)$}  \label{subsec.proofstep1} In this step, we describe the construction of $\nu (n,T)$.
We consider $n$ points $x_1, x_2, ...., x_n \in M$ 
and construct a closing-with curve $\nu (x_1, x_2, ...., x_n, T)$, i.e a curve such that $$\Gamma (x_1, x_2, ...., x_n, T) := \gamma (x_1, x_2, ...., x_n, T) \cup \nu (x_1, x_2, ...., x_n, T)$$ is a closed homologically trivial curve
with some uniform control of its length. Notice that at this stage, the construction will be holding for every choice of the points $x_1, x_2, ...., x_n \in M$ and that we will have to choose these points more carefully so that they satisfy all items of 
Proposition \ref{2-10*}.
We denote by $(\upsilon_i)_{1 \le i \le k}$ a basis of the torsion free part of $H_1(M,\ZZ)$ (not to be confused with $H_1(M, \ZZ)$) and by $(\beta_j)_{1 \le j \le k}$ a family of closed 1-forms dual to $(\upsilon_i)_{1 \le i \le k}$. Meaning that for any $ 1 \le i, j \le k$ we have
	$$ \inte{\upsilon_j} \beta_i = 1 \ \text{ if } \ i = j \text{ and } \ 0 \text{ otherwise }  \ . $$

The Step 1 can be formulated through the following Lemma.

\begin{lemma}\label{lemmafirststep}
	\label{lemmedetermcostclosingcurve}
There is a constant $C > 0$ such that	for any $\epsilon > 0$ there is $T > 0$ such that for any $n \in \NN$ and any $x_1, ... , x_n \in  M$ there is a curve $\nu(x_1, ... x_n,T)$ such that 
		\begin{enumerate}
			\item the curve
				 $$\Gamma(x_1, ... ,x_n, T) := \gamma(x_1, ... , x_n,T) \cup \nu(x_1, ... , x_n, T) $$ is closed and homologically trivial. 
			\item $$ \frac{ l(\nu(x_1, ... , x_n, T))}{nT} \le \epsilon + C \Big| \supr{1 \le j \le k} \, \frac{1}{nT} \inte{\gamma(x_1, ... , x_n , T)} \ \beta_j  \Big| \ .$$
		\end{enumerate}
\end{lemma}

\begin{figure}[h!]
\begin{center}
	\def\svgwidth{1 \columnwidth}
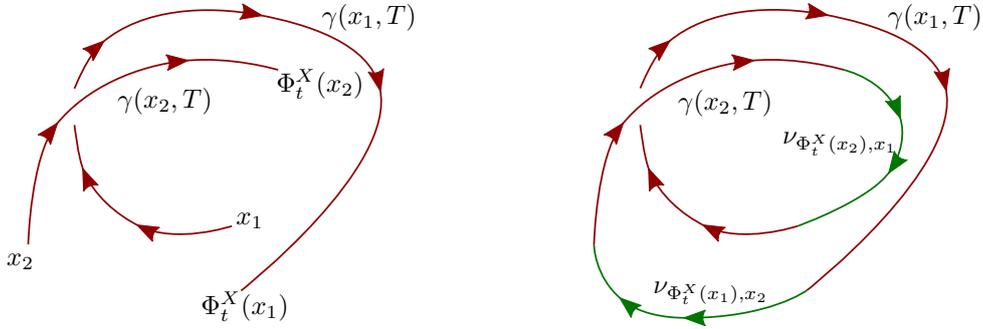
	\end{center}
\caption{On the left side, two trajectories of the $X$ flow up to time time $T$ starting at points $x_1$ and $x_2$. On the right, the closed up curve $\gamma_c(x_1, x_2,T)$.}
	\label{fig.closing2}
\end{figure}

\textbf{Proof.} The construction of $\nu(x_1, ... , x_n,T)$ is in  two steps.  First we glue altogether the pieces of $\gamma(x_1, ... , x_n, T)$ with curves of uniformly bounded length. In fact, for any $x, y \in M$, we choose a segment $\nu_{x,y}$ from $x$ to $y$ which is smoothly tangent to the foliation induced by $X$ at $x$ and $y$. One can moreover assume that 
	$$ \supr{x,y \in M} \ l(\nu_{x,y}) \le 2 D \ ,$$
where $D$ is the diameter of $M$. We set 
 $$ \tilde{\nu}(x_1, ... , x_n) := \underset{1 \le i \le n}{\cup} \nu_{\Phi_t^X(x_{i}),x_{i+1}} \ ,  $$
where the indices in the above union are to be understood with cyclic order: we glue $\gamma(x_n,T)$ with $\gamma(x_1, T)$ (see Figure \ref{fig.closing2}). Finally, we set 
	$$ \gamma_c(x_1, ... , x_n, T) = \gamma(x_1, ... , x_n , T) \cup \tilde{\nu}(x_1, ... , x_n) \ .$$

\begin{figure}[h!]
\begin{center}
	\def\svgwidth{0.8 \columnwidth}
%% Creator: Inkscape 1.0.1 (c497b03c, 2020-09-10), www.inkscape.org
%% Accompanies image file 'closingcurve4.pdf' (pdf, eps, ps)
\begingroup%
  \makeatletter%
  \providecommand\color[2][]{%
    \errmessage{(Inkscape) Color is used for the text in Inkscape, but the package 'color.sty' is not loaded}%
    \renewcommand\color[2][]{}%
  }%
  \providecommand\transparent[1]{%
    \errmessage{(Inkscape) Transparency is used (non-zero) for the text in Inkscape, but the package 'transparent.sty' is not loaded}%
    \renewcommand\transparent[1]{}%
  }%
  \providecommand\rotatebox[2]{#2}%
  \newcommand*\fsize{\dimexpr\f@size pt\relax}%
  \newcommand*\lineheight[1]{\fontsize{\fsize}{#1\fsize}\selectfont}%
  \ifx\svgwidth\undefined%
    \setlength{\unitlength}{1065.99545144bp}%
    \ifx\svgscale\undefined%
      \relax%
    \else%
      \setlength{\unitlength}{\unitlength * \real{\svgscale}}%
    \fi%
  \else%
    \setlength{\unitlength}{\svgwidth}%
  \fi%
  \global\let\svgwidth\undefined%
  \global\let\svgscale\undefined%
  \makeatother%
  \begin{picture}(1,0.4256568)%
    \lineheight{1}%
    \setlength\tabcolsep{0pt}%
    \put(0,0){\includegraphics[width=\unitlength,page=1]{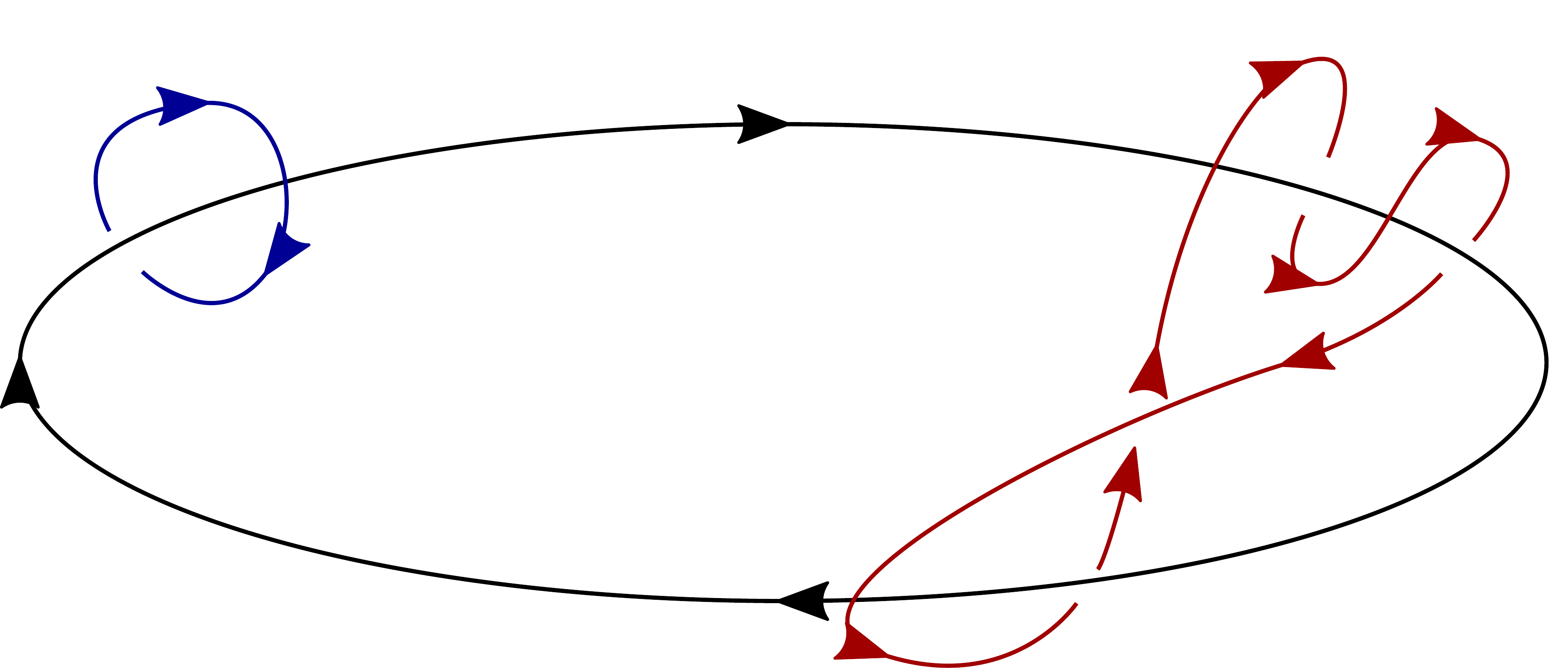}}%
    \put(0.0499742,0.38382681){\color[rgb]{0,0,0}\makebox(0,0)[lt]{\lineheight{1.25}\smash{\begin{tabular}[t]{l}$\nu$\end{tabular}}}}%
    \put(0.7552734,0.40372671){\color[rgb]{0,0,0}\makebox(0,0)[lt]{\lineheight{1.25}\smash{\begin{tabular}[t]{l}$\gamma_c(x_1,..., x_n,T)$\end{tabular}}}}%
    \put(0.38528757,0.37089189){\color[rgb]{0,0,0}\makebox(0,0)[lt]{\lineheight{1.25}\smash{\begin{tabular}[t]{l}$\mathbb{S}^1$\end{tabular}}}}%
    \put(0,0){\includegraphics[width=\unitlength,page=2]{closingcurve4.pdf}}%
  \end{picture}%
\endgroup%

	\end{center}
\caption{Here $M = \RR^3 \setminus \mathbb{S}^1$. The red curve is not homologicaly trivial since it links $-1$ times with $\mathbb{S}^1$. The blue curve  $\nu$ is any meridian of the knot defined by $\mathbb{S}^1$: its class generates $H_1(\RR^3 \setminus \mathbb{S}^1)$. The green dotted paths are meant to connect $\gamma(x_1,..., x_n,T)$ to the special representative $\nu$. The union of the (slitted) red curve, the two green paths, and the (slitted) blue curve travelled up as many times as the red one links the circle form the desired curve $\Gamma(x_1,..., x_n, T)$.}
	\label{fig.unrollhomology}
\end{figure}

Note that 
\begin{equation}
	\label{eq.closing1}
	l(\tilde \nu  (x_1,x_2,....,x_n, T)) =   l( \gamma_c(x_1, ... , x_n, T)) - l(\gamma(x_1, ... , x_n , T))  \le 2D n \ .
\end{equation}
The curve $\gamma_c(x_1, ... , x_n, T)$ has no reason to be homologically trivial. The second part of the proof is devoted to 'homologically unroll it'. In order to do so, we connect $\gamma_c(x_1, ... , x_n, T)$ to the previously chosen system of closed curves $(\upsilon_i)_{1 \le k \le k}$. Recall that it is a system of closed curves generating the ring of integer coefficients in $H_1(M, \RR)$. The closed curve $\gamma_c(x_1, ..., x_n, T)$, as a cycle with real coefficients, decomposes homologically with respect to the $(\upsilon_i)_{1 \le j \le k}$: there are integers $c_1, ... , c_k$ such that
	$$ \gamma_c(x_1, ..., x_n, T) - \somme{1 \le j \le k} \ c_j \upsilon_j $$ 
is homologically trivial (the above formula is to be understood homologically). We shall now make this homological relation geometric by relating the  $ \gamma_c(x_1, ..., x_n, T)$ with the $\upsilon_i$. We slit open the curve $ \gamma_c(x_1, ..., x_n, T)$ at some definite point and insert connecting paths from it to any of the $(\upsilon_i)_{1 \le i \le k}$ as in Figure \ref{fig.unrollhomology}. This can be done using a similar system of paths as defined to construct $\tilde{\nu}(x_1,..., x_n)$. As this part of the argument is completely similar we left the details to the reader. The resulting curve, the desired curve $\Gamma(x_1, ... , x_n, T)$, is then a real homologically trivial closed curve which writes as
	$$ \Gamma(x_1, ... , x_n, T) = \gamma(x_1, ..., x_n, T) \cup \nu(x_1, ... , x_n, T) \ ,$$ 
and satifies the item (1) of Lemma \ref{lemmafirststep}.\\

Note that the construction gives
\begin{equation}
\label{eqproofcurvehomlenghtdeffect}
 	\begin{split}
			 l(\nu(x_1 ,..., x_n, T)) & =  l(\Gamma(x_1, ... , x_n, T)) - l(\gamma(x_1, ..., x_n, T) \\
			  & \le 2  D n + 4 D k + \supr{ 1 \le j \le k }  l(\upsilon_j)  \cdot \supr{ 1 \le j \le k } | c_j | \ ,
	 \end{split}
\end{equation}
where 
\begin{itemize}
\item the term $2D n$, corresponding to the closing, comes from \eqref{eq.closing1};
\item the term $4 D k$ correspond to the length of the connecting system of paths coming back and forth from $\gamma(x_1, ..., x_n, T)$ to the $\upsilon_j$'s;
\item the term $\supr{ 1 \le j \le k }  l(\upsilon_j)  \cdot \supr{ 1 \le j \le k } | c_j | $ 	 corresponds to the length of the correcting homological curves $\underset{1 \le j \le k}{\cup} c_j \upsilon_j$. 
\end{itemize}

The final step consists in expressing the coefficients $(c_i)_{1 \le i \le k}$ in term of the closed $1$-forms $(\beta_j)_{1 \le j \le k}$ whose classes generate $H^1_{\mathrm{DR}}(M, \RR)$. Together with 
(\ref{eqproofcurvehomlenghtdeffect}), this will prove item (2) of Lemma \ref{lemmafirststep}.
Let us show that there is a constant $C_2$ such that for any $1 \le j \le k $ one has
	$$| c_j | \le C_2 \left( 2 D n + \Big|\inte{\gamma(x_1, ... , x_n , T)} \ \beta_j \Big| \right) \ .  $$ 

The proof is straightforward and consists in integrating the forms $\beta_j$ over the homologically trivial curve 
	$$\widetilde{\Gamma} := \gamma_c(x_1, ... ,x_n, T) \underset{1 \le j \le k}{\cup} c_j \upsilon_j \ .$$

On the one hand one has for all $1 \le j \le k$ 
	$$ \inte{\widetilde{\Gamma}} \beta_j = 0 \ . $$
the above integral also writes as 
\begin{align*}
	\inte{\widetilde{\Gamma}} \beta_j & =  \inte{\gamma_c(x_1, ... ,x_n, T)} \beta_j  +  \somme{1 \le i \le k} \ c_ i \inte{\upsilon_i} \beta_j  \\
		& = \inte{\gamma_c(x_1, ... ,x_n, T)} \beta_j  +  c_j \ ,
\end{align*}
by construction of $\beta_j$. Therefore,
\begin{align*}
	|c_j| & = \Big| \inte{\gamma_c(x_1, ... ,x_n, T)} \beta_j \Big| \le  \Big| \inte{\gamma(x_1, ... ,x_n, T)}  \beta_j \Big| + \Big| \inte{\tilde{\nu}(x_1, ... ,x_n, T)}  \beta_j \Big| \\
		 & \le \Big| \inte{\gamma(x_1, ... ,x_n, T)}  \beta_j \Big| + C_3 \ l(\tilde{\nu}(x_1, ..., x_n, T)) \\ 
		 & \le  \Big| \inte{\gamma(x_1, ... ,x_n, T)}  \beta_j \Big| + C_3 2 D n \ ,
\end{align*}
where 
	$$C_3 := \supr{1 \le j \le k} ||\beta_j||_{\infty}  \ . $$

Looking backward to \eqref{eqproofcurvehomlenghtdeffect} we get 

 \begin{align*}
 	  l(\nu(x_1 ,..., x_n, T))   & \le 2  D n + 4 D k + \supr{ 1 \le j \le k }  l(\upsilon_j)  \cdot \supr{ 1 \le j \le k } | c_j | \\ 
 	  & \le  2  D n + 4 D k + \supr{ 1 \le j \le k }  l(\upsilon_j) \supr{1 \le j \le k} \left( \Big| \inte{\gamma(x_1, ... ,x_n, T)}  \beta_j \Big| + C_3 2 D n \  \right) \\
 	  & \le C_4 n  + C_5 \supr{1 \le j \le k} \Big| \inte{\gamma(x_1, ... ,x_n, T)}  \beta_j \Big| \ ,
\end{align*} 
for some constants $C_4$ and $C_5$. Therefore one has 
	$$ \frac{  l(\nu(x_1 ,..., x_n, T)) }{nT} \le  \frac{C_4}{T} + C_5\supr{1 \le j \le k} \Big| \frac{1}{nT} \inte{  \gamma(x_1, ... ,x_n, T)} \ \beta_j \Big| \ . $$
Setting $T$ such that $C_4/T \le \epsilon$ concludes.  \hfill $\blacksquare$ \\

\subsection{Step 2: the starting points} \label{subsec.proofstep2}

Let us then come back to the proof of Proposition \ref{2-10*}. The idea is to consider the set of pieces of trajectories $\gamma (x_1,x_2,....,x_n, T)$ 
of $X$ as a probability space and show the positivity of the probability of the
set of trajectories such that the closed curve $\Gamma (n,T)$ constructed in step 1 satisfy all items of the Proposition \ref{2-10*}. In order to do so we fix $\epsilon > 0$ and $T$ such that the conclusions of Lemma \ref{lemmedetermcostclosingcurve} hold. We set
	$$ (\Omega, \mathcal{A}, \PP) := (M^{\NN}, \mathcal{A}, \mu^{\NN}) \ ,$$ where $\mathcal{A}$ is the cylindrical $\sigma$-algebra on the product space $M^{\NN}$ and  $\mu^{\NN}$ is the product measure.
 We denote by $\omega$ an event of $\Omega$, by $Z_i$ the projection that associates to an event $\omega$ its $(i-1)$th factor. The $(Z_i)_{i > 0}$ are by construction independent identically distributed random variables taking value in $M$: they correspond to the random starting points. For any $T > 0$ we denote by 
	$$\gamma(n,T, \omega) = \gamma(Z_1(\omega), ... , Z_n(\omega),T) \ . $$

The following simple lemma is the key. Its proof is a basic application of the (weak) law of large number.
\begin{lemma}
	\label{lemmaproofproprayleyL1LNN}
Let $M$ be a closed Riemannian manifold. Let $\epsilon, T > 0$ and $\eta$ be a bounded 1-form. Then the following holds
		$$ \PP \left( \Big|  \frac{1}{nT}  \inte{\gamma(n, T, \,.)} \eta  - \inte{M} \ i_X(\eta) \ d \mu \Big| \le \epsilon \right) \tends{ n \to \infty} 1 \ . $$
\end{lemma}

\textbf{Proof.} For any fixed differential form $\eta$ and any $T > 0$ we define
	$$U_i^{\eta, T}(\omega) := \frac{1}{T} \inte{\gamma(Z_i(\omega),T)} \eta \ . $$
The above defined random variables are independent and identically distributed random variables since they write $f(Z_i)$ for some measurable function $f$, where the $(Z_i)_{1 \le i \le n}$ are themselves independant 
identically distributed random variables by definition. 
Notice also that the random variables $U_i^{\eta, T}$ are bounded since $\eta$ is assumed to be bounded and since the curves $\gamma(x,T)$ have uniformly in $x$ bounded length. \\

Note that $$\frac{1}{nT}  \inte{\gamma(n, T, \omega)} \eta = \frac{1}{n} \somme{1 \le i \le n} U_i^{\eta, T}(\omega) $$ by construction of $\gamma(n,T, \omega)$. Therefore 
$$ \frac{1}{nT}  \inte{\gamma(n, T, \omega)} \eta $$ writes as the mean of independant identically distributed bounded random variables. 
One can then apply the law of large number to get the conclusion of Lemma \ref{lemmaproofproprayleyL1LNN} provided that 
	$$ \EE( U_i^{\eta, T} )= \inte{M} \ i_X \eta \ d \mu \ .$$
The above equality follows from a simple computation. Indeed, for any $i \ge 1$ one has 
\begin{align*}
	 \EE \left( \frac{1}{T} \inte{\gamma(Z_i, T)} \eta \right) & = \inte{M} \frac{1}{T} \inte{\gamma(x, T)} \eta \ d \mu(x) \\
	 & = \inte{M} \frac{1}{T} \inte{0}^T \eta(X(\Phi_t^X(x)) \ dt \,d \mu(x) \\
	 & = \inte{M} \frac{1}{T} \inte{0}^T (i_X\eta)(\Phi_t^X(x)) \ dt \, d \mu(x) \\
	 & = \frac{1}{T}\inte{0}^T \inte{M}  (i_X\eta)(\Phi_t^X(x)) \ d \mu(x) \, dt\\
	 & = \inte{M}  (i_X\eta) \ d \mu \ , 
\end{align*}
since $\mu$ is invariant under $\Phi_t^X$ by construction of $X$. \hfill $\blacksquare$ \\

Note that we have a stronger convergence than the one stated in Lemma \ref{lemmaproofproprayleyL1LNN} since the (strong) law of large number asserts that an almost sure convergence holds. Note also that we could have stated Lemma \ref{lemmaproofproprayleyL1LNN} with the assumption that $\eta$ is integrable but we would also have had to show that the random variables $U_i^{\beta, T}$ are integrable. \\

We are now ready to find points $x_1, x_1, ...,x_n$ such that for every $T>0$, the curve $\Gamma (n, T)$ constructed in step 1 and satisfying item $3$, also satisfies items $1$, $2$ and $4$ of Proposition \ref{2-10*}. In order to do it, 
we shall apply Lemma \ref{lemmaproofproprayleyL1LNN} to several 1-forms according to each item. \\

{\bf Item 1:}
 Recall that $(\beta_1, \beta_2, ... , \beta_k)$ have been defined as a smooth 1-forms basis of $H^1_{\mathrm{DR}}(M, \RR)$. Lemma \ref{lemme.interiorproduct} implies that for all $ 1 \le j \le k$
$$	\inte{M} \ i_X \beta_j = \inte{M} d \alpha \wedge \beta_j = 0 \ , $$
since $\beta_j$ is closed. Lemma \ref{lemmaproofproprayleyL1LNN} then gives that for every $1 \le j \le k$ one has 
	\begin{equation}
			\label{eqlemproofLNN1}
			 \PP \left( \Big| \frac{1}{nT}  \inte{\gamma(n, T, \,.)} \beta_j \Big| \le \epsilon \right) \tends{ n \to \infty} 1 \ . 
	\end{equation}
	
We set 
		$$ \nu(n,T,\omega) := \nu(Z_1(\omega), ... , Z_n(\omega),T) \ ,$$
where $\nu$ is the closing-with curve given by Lemma \ref{lemmedetermcostclosingcurve}. \\

Note that if $\omega \in \Omega$ is such that for every $1 \le j \le k$ we have
$$ \Big| \frac{1}{nT}  \inte{\gamma(n, T, \omega)} \beta_j \Big| \le \epsilon $$
then Lemma \ref{lemmedetermcostclosingcurve} (2) implies
\begin{equation}
			\label{eqlemproofLNN1*}
 \frac{ l(\nu(n, T, \omega))}{nT} \le (1 + C)  \epsilon \ .
\end{equation}

From (\ref{eqlemproofLNN1}) and (\ref{eqlemproofLNN1*}) we get
\begin{equation}
			\label{eqlemproofLNN1**}
			 \PP \left( \frac{ l(\nu(n, T, \, .))}{nT} \le (1 + C)  \epsilon \right) \tends{ n \to \infty} 1 \ 
	\end{equation}

which implies that, for $n$ large enough and any $\omega$ in a subset of $\Omega$ of probability close to $1$, the first item of Proposition \ref{2-10*} is satisfied.\\

{\bf Item 2:} 
We now apply Lemma \ref{lemmaproofproprayleyL1LNN} to the differential form 
	 $$\widetilde{ \star d \alpha}  := 
		\left\{ 
			\begin{array}{l}
				\frac{ (\star d \alpha)_x}{| (\star d \alpha)_x|} \ \text{ if } \ (\star d \alpha)_x \neq 0 \\
				0 \ \text{otherwise} \ .
			\end{array} 
		\right. $$
		
Lemma \ref{lemme.interiorproduct} implies that 
\begin{equation}\label{eqitem2-1}
 \inte{M} i_X\left( \widetilde{ \star d \alpha} \right)  = \inte{M} \widetilde{ \star d \alpha}  \wedge d \alpha = \inte{M} \ |\star d \alpha| \ d \mu = || d \alpha||_{1},
\end{equation}
the last two equality coming from the definition of $\widetilde{ \star d \alpha}$.
%On the one hand, the integral appearing in Lemma \ref{lemmaproofproprayleyL1LNN} gives 
		%$$ \inte{M} \widetilde{ \star d \alpha}  \wedge d \alpha = \inte{M} \ |\star d \alpha| \ d \mu = || d \alpha||_{1} \ .$$
On the other hand, the integral over the random curves appearing in Lemma \ref{lemmaproofproprayleyL1LNN} in the summation values of the term $$\inte{ \gamma(n,T, \, .)} \widetilde{ \star d \alpha}$$ 
writes
$$\inte{ \gamma(Z_i (.),T)} \widetilde{ \star d \alpha} = \inte{0}^T \ |\widetilde{ \star d \alpha}|(\Phi_t^X(Z_i (.))) \ dt = l(\gamma(Z_i (.),T)) \ . $$
and therefore Lemma \ref{lemmaproofproprayleyL1LNN} implies 
	\begin{equation}
			\label{eqlemproofLNN2*}
			 \PP \left( \Big| \frac{l(\gamma(n,T,\cdot)}{nT} - || d \alpha||_{1}  \Big| \le \epsilon \right) \tends{ n \to \infty} 1 \ . 
	\end{equation}
 As above, we conclude that for $n$ large enough and any $\omega$ in a subset of $\Omega$ of probability close to $1$, the item $2$ of Proposition \ref{2-10*} is satisfied.\\

%The event above corresponds to item (2) of Proposition \ref{propexistencecurve}. \\
{\bf Item 4:} 	
We finally apply Lemma \ref{lemmaproofproprayleyL1LNN} to a $1$-form $\beta$ such that $d \beta = \star \alpha$.
%which corresponds to item (4) of Proposition \ref{propexistencecurve}. 
We use Stokes formula to get 
	$$ \inte{M} \beta \wedge d \alpha = \inte{M} d \beta \wedge \alpha = \inte{M} \star \alpha \wedge \alpha = || \alpha||^2_{2} \ .  $$
Lemma \ref{lemmaproofproprayleyL1LNN} then gives 
\begin{equation}
			\label{eqlemproofLNN3}
			 \PP \left( \Big| \frac{1}{nT}  \inte{\gamma(n, T, \omega)} \beta - || \alpha||^2_{2} \ d \mu \Big| \le \epsilon \right) \tends{ n \to \infty} 1 \ . 
\end{equation} 
We conclude that for $n$ large enough and any $\omega$ in a subset of $\Omega$ of probability close to $1$, the item $4$ of Proposition \ref{2-10*} is satisfied.\\

To summarize, by \eqref{eqlemproofLNN1**} \eqref{eqlemproofLNN2*} and \eqref{eqlemproofLNN3} each item $1$, $2$ and $4$ corresponds to an event 
%by \eqref{eqlemproofLNN1} \eqref{eqlemproofLNN2} and \eqref{eqlemproofLNN3} 
having probability to be satisfied which tends to $1$ as $n \to \infty$. An element in the intersection of all these events satisfies all the items of Proposition \ref{2-10*}. 
%Note that there is only finitely many events under consideration (as many as the dimension of $H^1_{\mathrm{DR}}(M, \RR)$ for \eqref{eqlemproofLNN1} + 2 corresponding to \eqref{eqlemproofLNN2} and \eqref{eqlemproofLNN3}). 
Therefore, for $n$ large enough, all these events must intersect non trivially. \\

In particular, there exists $n \in \NN^*$ and finitely many starting points $x_1, ... , x_n \in M$ such that the following union of curves
$$ \gamma(n,T) := \underset{1 \le i \le n}{\cup} \ \gamma(x_i,T) $$ 
together with its closing curve $\nu(x_1, ..., x_n,T)$ satisfies all the items of Proposition \ref{2-10*}, concluding. \hfill $\blacksquare$

\section{An example: the collapsing Hopf fibration}\label{section4}
\label{sec.example}

Consider the $\mathbb{S}^1$-isometric action on the round $3$-sphere $$\mathbb{S}^3 =\{ (z_1,z_2) \in \mathbb{C}^2\,|\, |z_1|^2 +|z_2|^2 =1 \}$$ given by 
\begin{equation}\label{hopf}
e^{i\theta} (z_1, z_2) := (e^{i\theta} z_1, e^{i\theta}z_2).
\end{equation}
This action provides $\mathbb{S}^3$ with a $\mathbb{S}^1$-fiber bundle structure $\mathbb{S}^1 \rightarrow \mathbb{S}^3 \rightarrow \mathbb{CP}^1$ over the one dimensional complex projective space $\mathbb{CP}^1$ called the Hopf fibration. The projection $\pi : \mathbb{S}^3 \rightarrow \mathbb{CP}^1$
turns out to be a Riemannian submersion over $\mathbb{CP}^1$ endowed with the round metric of constant sectional curvature equal to $4$. We now endow $\mathbb{S}^3$ with the family of Berger's metrics $(g_{\epsilon})_{0 < \epsilon \le 1}$ defined as deformations along the Hopf fibers of the round metric $g_1$ as follows. The tangent bundle $T \mathbb{S}^3 =V \oplus H$ splits as the sum of the one dimensional 'vertical' subbundle $V$ tangent to the Hopf fibers and its 'horizontal' $g_1$-orthogonal complement $H$. This gives rise to the following decomposition of the metric $g_1 = g_v + g_h$, where $g_v$ (resp. $g_h$) denotes the restriction of the metric $g_1$ to $V$ (resp. to $H$). The family of Berger's metrics $g_{\epsilon}$ is the following deformation of $g_1$,
	$$ g_{\epsilon} := g_h + \epsilon^2 g_v \ . $$
We emphasize with a subscript the different quantities associated to the $g_{\epsilon}$ metric (for example, the length of a curve $\gamma$ for the metric $g_{\epsilon}$ will be denoted $l_{\epsilon}(\gamma)$). Endowed with the Gromov-Hausdorff topology, the family of Berger spheres $(S^3, g_\epsilon)$ converges  to the 2-sphere endowed with the round metric. This family exhibits three main features:
\begin{itemize}
\item it collapses as the injectivity radius at any point tends to $0$ as $\epsilon \to 0$;
\item it has uniformly bounded from above (and below) diameter (it actually Gromov-Haussdorff converges  toward $\mathbb C \mathbb P ^1$);
\item it has uniformly bounded sectional curvature.
\end{itemize} 
In particular, the metrics $g_{\epsilon}$ satisfy to all the assumptions of Theorem \ref{maintheo}. The $\mathbb{S}^1$-action (\ref{hopf}) is isometric for each of the Berger's metrics $g_\epsilon$ hence its generating vector field $X$ is divergence free. Let us consider $\alpha := X^\flat$ 
the dual $1$-form of $X$ with respect to $g_1$. Notice that $\epsilon ^{2} \alpha$ is the dual $1$-form of $X$ for the metric $g_\epsilon$ and that $\alpha$ is therefore a coexact $1$-form with respect to every metric $g_\epsilon$. One can verify that $\Delta _1 \alpha = 4\epsilon ^2 \alpha$ as well as that $$\lambda ^1 _0 (\epsilon):= 4\epsilon ^2$$ is the smallest possible eigenvalue of the $1$-forms on $(S^3, g_\epsilon)$. Another important property of $\alpha$ is that 
\begin{equation}\label{d-alpha}
d\alpha = 2\pi ^* \Omega,
\end{equation}
where $\Omega$ is the volume form of $\mathbb{CP}^1$, that is the $2$-form on  $\mathbb{CP}^1$ satisfying $$\int_{ \mathbb{CP}^1} \Omega =  \vol  (\mathbb{CP}^1) = \pi.$$

The goal of this subsection is to show that Theorem \ref{theo.cheegersurfaceminimale} is sharp along the collapsing in the sense of the following
\begin{proposition}
\label{asympto-h}
With the notations introduced above, there is a constant $a >0$ such that we have
	$$  \lim _{\epsilon \to 0} \frac{h^1_{\epsilon}}{\epsilon} = a  \ . $$
\end{proposition}

\textbf{Proof.} We start by showing an upper bound on $h_{\epsilon}$. This upper bound will be improved latter into the asymptotic appearing in the conclusion of Proposition \ref{asympto-h}.

\begin{lemma}
	\label{lemme.upper}
With the notation introduced above, we have for all $\epsilon \in ]0,1]$
$$  h^1_\epsilon \leq 2\epsilon \ .$$
\end{lemma}

\textbf{Proof.} Recall that 
	$$ h^1_{\epsilon} := \infi{\gamma \hookrightarrow M} \ \frac{l_{\epsilon}(\gamma)}{A_{\epsilon}(\gamma)} \ , $$

therefore, in order to show the upper bound, it is sufficient to find a curve with the above isoperimetric ratio bounded above. The curve we will consider is any given fiber of the Hopf fibration that we denote simply by $\gamma$. 
By definition of the metric $g_\epsilon$ its length $l_{\epsilon} (\gamma)$ is easy to compute:
	$$ l_{\epsilon}(\gamma) = 2 \pi \epsilon \ .$$
Therefore we have 
\begin{equation}\label{upper1}
h^1_{\epsilon} \le \frac{2 \pi \epsilon}{A_{\epsilon}(\gamma)} \ ,
\end{equation}
hence, in order to prove \ref{lemme.upper}, one is left to show that $A_{\epsilon}(\gamma) \geq \pi$. This will be a consequence of the following observation. Let $S$ be a surface whose boundary is $\gamma$. Since  $\pi : \mathbb{S}^3 \rightarrow \mathbb{CP}^1$ 
is a Riemannian submersion, we have for any $x \in \mathbb{S}^3$ and any $\epsilon > 0$ 
	$$| \pi ^* \Omega |_{\epsilon} (x) =1 \ ,$$
and then 
$$
\int_{S} \pi ^* \Omega \leq  \area_{\epsilon}( S) \ .
$$
We then deduce by (\ref{d-alpha}) and Stokes formula, 
\begin{equation*}
\frac{1}{2} \int _\gamma \alpha = \int_{S} \pi ^* \Omega \leq  \area_{\epsilon} (S) \ ,
\end{equation*}
thus, 
\begin{equation*}
\pi \leq \area_{\epsilon} (S) \ ,
\end{equation*}
since by definition of $\alpha$, we have 
$$\int_\gamma \alpha = 2\pi \ .$$

Now, $S$ being arbitrary in the above inequality, we get 
\begin{equation}\label{upper2}
\pi \leq  A_\epsilon (\gamma) \ .
\end{equation}
Therefore, using \eqref{upper1}, we have 
$$
h^1_\epsilon \le \frac{2 \pi \epsilon}{A_{\epsilon}(\gamma)} \le 2\epsilon \ ,
$$
which concludes the proof of Lemma \ref{lemme.upper}. \hfill $\blacksquare$ \\

The next Lemma shows that a family of curves $\gamma_\epsilon$ which almost realize the constant $h^1_{\epsilon}$ is 'almost vertical'; we say that a family of curve $(\gamma_{\epsilon})_{\epsilon \in ]0,1]}$ is \textbf{almost realizing} $h^1_{\epsilon}$ 
if for all $\epsilon > 0$ we have
	$$ \frac{l_{\epsilon} (\gamma_{\epsilon})}{A_{\epsilon}(\gamma_{\epsilon})} \le  h^1_{\epsilon} (1+ \epsilon) \ . $$
Given a curve $\gamma : [0,1] \to S^3$, the tangent vector $\dot{\gamma}(t) = \dot{\gamma} ^h(t)  + \dot{\gamma}^v (t)  $ decomposes along
the vertical and horizontal bundles $V$ and $H$ and we call the horizontal (resp. vertical) length of $\gamma$ with respect to the metric $g_\kappa$, $\kappa \in (0,1[$ the quantity
$$l_\kappa^h (\gamma) := \int_0 ^1 |\dot{\gamma}^h (t) | _\kappa dt , \, \, \, \, \, \mbox{resp.} \, \, \,\,\,\,  l_\kappa^v (\gamma) := \int_0 ^1 |\dot{\gamma}^v (t) | _\kappa dt. $$

\begin{lemma}
\label{vertical}
For any $\kappa$ there is a constant $C_{\kappa}$ such that for any family $\gamma_{\epsilon}$ almost realizing $h^1_{\epsilon}$ with $0<\epsilon \leq \kappa$, we have
	$$ l_{\kappa}(\gamma_{\epsilon}) \le  (1+ C_{\kappa} \cdot \epsilon )\,\, l^v_\kappa(\gamma_{\epsilon})  .  $$
\end{lemma}

\textbf{Proof.}  By the triangle inequality we have,
	$$ l_{\kappa}(\gamma_{\epsilon}) \le  l^v_{\kappa}(\gamma_{\epsilon}) + l^h_{\kappa}(\gamma_{\epsilon}) \ , $$

so that we will conclude provided that there is a constant $C_\kappa$ such that for all $\epsilon \leq \kappa$,
 \begin{equation}\label{equiv-h-0}
 l^h_{\kappa}(\gamma_{\epsilon}) \le C_\kappa \  \epsilon \ l_\kappa(\gamma_{\epsilon}) \ .
 \end{equation}
 We will now focus on proving (\ref{equiv-h-0}). We start with the observation that for any $\kappa, \epsilon \in ]0,1]$ we have
$$ l^h_{\kappa}(\gamma_{\epsilon}) = l^h_{\epsilon}(\gamma_{\epsilon}) \le l_{\epsilon}(\gamma_{\epsilon}) \ . $$
Since we assumed that $\gamma_{\epsilon}$ is almost realizing $h^1_{\epsilon}$ we also have 
	$$  l_{\epsilon}(\gamma_{\epsilon}) \le  \ h^1_{\epsilon} \ A_{\epsilon}(\gamma_{\epsilon}) \, (1+\epsilon) \ , $$
which implies, together with the already obtained upper bound \ref{lemme.upper}
	%$$ l_{\epsilon}(\gamma_{\epsilon}) \le 2 \ \epsilon  \ A_{\epsilon}(\gamma_{\epsilon}) \,(1+\epsilon) \ , $$
	%therefore we deduce 
\begin{equation}\label{equiv-h-1}
l_{\kappa}^h(\gamma_{\epsilon}) \le 2 \ \epsilon  \ A_{\epsilon}(\gamma_{\epsilon}) \,(1+\epsilon) \ .
\end{equation}
We denote by $S_{\gamma _ \epsilon}$ a minimal surface realising $A_{\epsilon}(\gamma_{\epsilon})$, i.e 
$$
A_{\epsilon}(\gamma_{\epsilon}) = \area_{\epsilon} (S_{\gamma _\epsilon}) \ . $$
Because $g_{\epsilon} \le g_\kappa$ for every $\epsilon \leq \kappa$, we have
	$$ \area_{\epsilon}(S_{\gamma_ \epsilon}) \le \area_{\kappa}(S_{\gamma _\epsilon}) \ ,$$
and, by definition of $h^1_{\kappa}$, we also have 
$$ \area_{\kappa}(S_{\gamma_\epsilon}) \le (h^1_{\kappa})^{-1} \cdot l_\kappa( \gamma_{\epsilon}) \ ,$$
therefore we obtain
	\begin{equation}\label{equiv-h-2}
	A_{\epsilon}(\gamma_{\epsilon}) \le (h^1_{\kappa})^{-1} \cdot l_\kappa( \gamma_{\epsilon}) \ .
	\end{equation}
Putting together the relations (\ref{equiv-h-1}) and (\ref{equiv-h-2}) yields 
	$$ l^h_{\kappa}(\gamma_{\epsilon}) \le 2 \ \epsilon \ (h^1_{\kappa})^{-1} \cdot l_\kappa( \gamma_{\epsilon}) \,(1+\epsilon) \, $$
	
which concludes the proof of (\ref{equiv-h-0}) setting $$C_\kappa := 4 \cdot (h^1_{\kappa})^{-1} \ge 2 \cdot (1+\epsilon) \cdot (h^1_{\kappa})^{-1} \ . $$ This ends the proof of Lemma \ref{vertical}. \hfill $\blacksquare$  \\

We are now ready to conclude the proof of Proposition \ref{asympto-h}. \\

We consider a family of curves $\gamma _\epsilon$ almost realizing $h^1_{\epsilon}$ as defined above. For $\epsilon \leq \kappa$, we also choose surfaces 
$S^\kappa_{\gamma_{\epsilon}}$ whose boundary is $\gamma_{\epsilon}$ and which realizes $A_\kappa (\gamma _\epsilon)$, i.e such that
$$ \area_{\kappa} (S^\kappa _{\gamma _\epsilon})  = A_\kappa (\gamma _\epsilon) \ . $$
By definition we have 
$$ h^1_{\kappa} \le \frac{l_\kappa(\gamma_{\epsilon}) }{ A_\kappa({\gamma_{\epsilon}})} = \frac{l_\kappa(\gamma_{\epsilon}) }{ \area_{\kappa} (S^\kappa  _{\gamma_{\epsilon}})} \ . $$ 
Using successively the above Lemma \ref{vertical} and the inequality $g_{\epsilon} \le g_\kappa$ we get 
$$ h^1_{\kappa}  \le  \frac{(1+ C_\kappa \cdot \epsilon) \,  \cdot l^v_\kappa(\gamma_{\epsilon}) }{ \area_{\kappa}(S^\kappa_{\gamma_{\epsilon}})}  \le  \frac{(1+ C_\kappa \cdot \epsilon) \cdot l^v_\kappa(\gamma_{\epsilon}) }{ \area_{\epsilon}(S^\kappa_{\gamma_{\epsilon}})} . $$                                 

But since $l^v_\kappa(\gamma_{\epsilon}) = \frac{\kappa}{\epsilon} l^v_{\epsilon}(\gamma_{\epsilon})$ we get 
$$ h^1_{\kappa}  \le  \frac{\kappa}{\epsilon} \frac{l^v_{\epsilon}(\gamma_{\epsilon})} { \area_{\epsilon}(S^\kappa_{\gamma_{\epsilon}})}\, (1+ C_\kappa \cdot \epsilon) \ . $$ 
Therefore, since $l_\epsilon ^v (\gamma _\epsilon) \leq l_\epsilon (\gamma _\epsilon)$ and $ \area_{\epsilon}(S^\kappa_{\gamma_{\epsilon}}) \geq A_\epsilon (\gamma _\epsilon)$,
$$ h^1_{\kappa}  \le \frac{\kappa }{\epsilon} \cdot  \frac{l_{\epsilon}(\gamma_{\epsilon})}{ A_{\epsilon}(\gamma_{\epsilon})} \, (1+ C_\kappa \cdot \epsilon) \ . $$ 
Recall that we have chosen $\gamma _\epsilon$ almost realizing $h^1_{\epsilon}$, thus we deduce from the above inequality
$$ h^1_{\kappa} \le \frac{\kappa}{\epsilon} \cdot h^1_{\epsilon} \,(1+\epsilon) \, (1+ C_\kappa \cdot \epsilon) \ .$$
We then have, for every $0<\epsilon \leq \kappa \leq 1$,
\begin{equation}\label{monotone}
\frac{ h^1_{\kappa}}{\kappa} \le \frac{ h^1_{\epsilon}}{\epsilon} \,(1+ \eta _\kappa (\epsilon)),
\end{equation} 
where $\lim_{\epsilon \to 0} \eta _\kappa (\epsilon) =0$.
We end the proof of Proposition \ref{asympto-h} by noticing that the above almost monotonicity property (\ref{monotone}) implies that 
$ \limsup_{\epsilon \to 0}  \frac{ h^1_{\epsilon}}{\epsilon} = \liminf_{\epsilon \to 0}  \frac{ h^1_{\epsilon}}{\epsilon} $ and that 
$\lim_{\epsilon \to 0}  \frac{ h^1_{\epsilon}}{\epsilon}  =a\ge h^1_1 >0$. \hfill $\blacksquare$

\section{A cusp-like counterexample}\label{section5}
\label{sec.counterexample}
This subsection is devoted to exhibit a counter example to Theorem \ref{maintheo} if one removes either the assumption that $M$ has bounded diameter or that $M$ has a Ricci curvature lower bound. We shall actually focus on exhibiting a counter example when $M$ has Ricci curvature bounded below but with a large diameter. Since our counterexample will carry negative curvature it also shows, by scaling the metric, that the lower bound on the Ricci curvature is mandatory. \\

\begin{figure}[h!]
\begin{center}
	\def\svgwidth{0.8 \columnwidth}
%% Creator: Inkscape 1.0.1 (c497b03c, 2020-09-10), www.inkscape.org
%% Accompanies image file 'tubes1.pdf' (pdf, eps, ps)
\begingroup%
  \makeatletter%
  \providecommand\color[2][]{%
    \errmessage{(Inkscape) Color is used for the text in Inkscape, but the package 'color.sty' is not loaded}%
    \renewcommand\color[2][]{}%
  }%
  \providecommand\transparent[1]{%
    \errmessage{(Inkscape) Transparency is used (non-zero) for the text in Inkscape, but the package 'transparent.sty' is not loaded}%
    \renewcommand\transparent[1]{}%
  }%
  \providecommand\rotatebox[2]{#2}%
  \newcommand*\fsize{\dimexpr\f@size pt\relax}%
  \newcommand*\lineheight[1]{\fontsize{\fsize}{#1\fsize}\selectfont}%
  \ifx\svgwidth\undefined%
    \setlength{\unitlength}{1275.59055118bp}%
    \ifx\svgscale\undefined%
      \relax%
    \else%
      \setlength{\unitlength}{\unitlength * \real{\svgscale}}%
    \fi%
  \else%
    \setlength{\unitlength}{\svgwidth}%
  \fi%
  \global\let\svgwidth\undefined%
  \global\let\svgscale\undefined%
  \makeatother%
  \begin{picture}(1,0.46666667)%
    \lineheight{1}%
    \setlength\tabcolsep{0pt}%
    \put(0,0){\includegraphics[width=\unitlength,page=1]{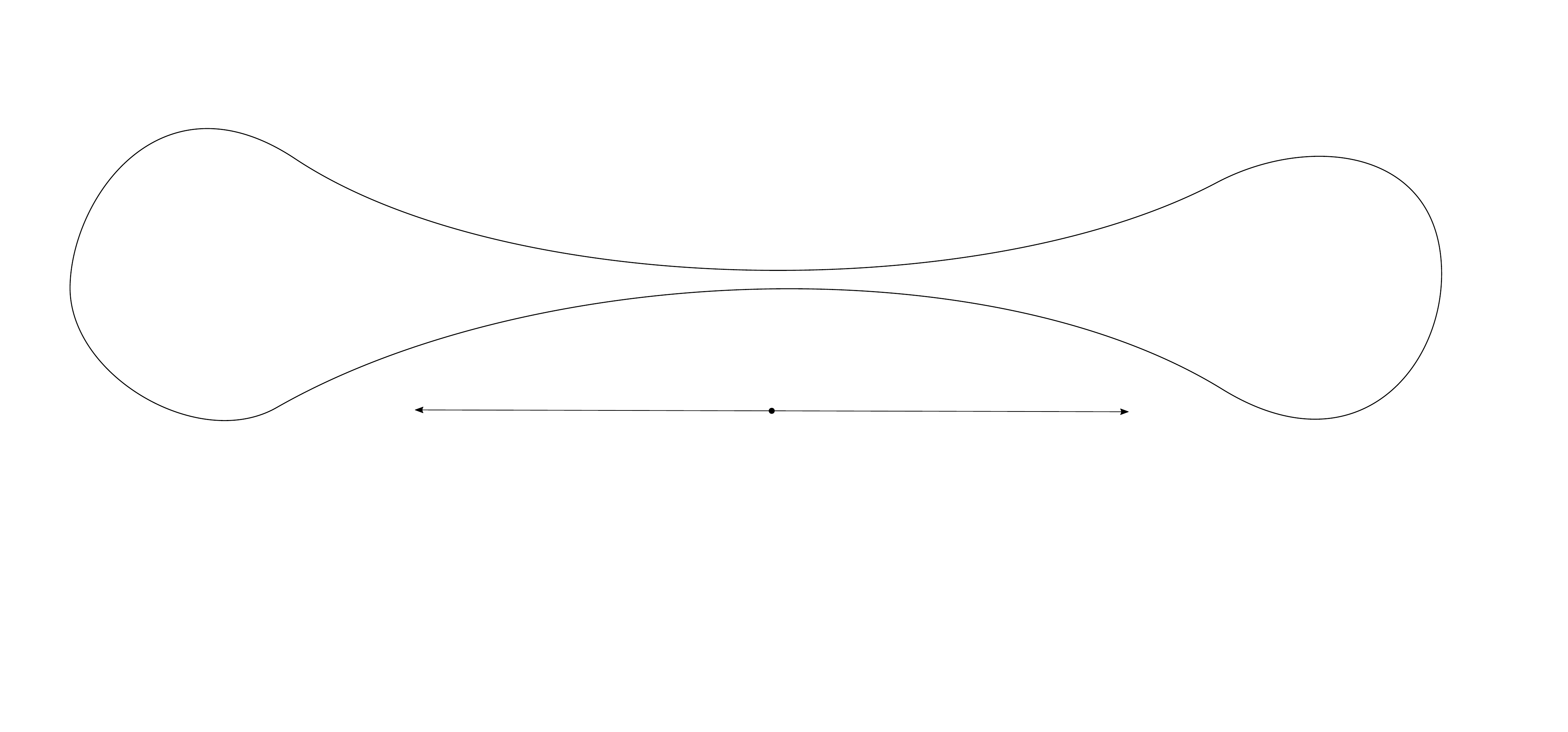}}%
    \put(0.47037035,0.17503703){\color[rgb]{0,0,0}\makebox(0,0)[lt]{\lineheight{1.25}\smash{\begin{tabular}[t]{l}$0$\end{tabular}}}}%
    \put(0.6923204,0.1679815){\color[rgb]{0,0,0}\makebox(0,0)[lt]{\lineheight{1.25}\smash{\begin{tabular}[t]{l}$\ln(1/\epsilon)$\end{tabular}}}}%
    \put(0,0){\includegraphics[width=\unitlength,page=2]{tubes1.pdf}}%
    \put(0.05056481,0.16327776){\color[rgb]{0.54901961,0.15686275,0}\makebox(0,0)[lt]{\lineheight{1.25}\smash{\begin{tabular}[t]{l}$M_L$\end{tabular}}}}%
    \put(0.84196296,0.16327779){\color[rgb]{0.4745098,0,0}\makebox(0,0)[lt]{\lineheight{1.25}\smash{\begin{tabular}[t]{l}$M_R$\end{tabular}}}}%
    \put(0.44006307,0.31744227){\color[rgb]{0.54901961,0.38039216,0}\makebox(0,0)[lt]{\lineheight{1.25}\smash{\begin{tabular}[t]{l}$M_C$\end{tabular}}}}%
    \put(0,0){\includegraphics[width=\unitlength,page=3]{tubes1.pdf}}%
    \put(0.2434167,0.17033335){\color[rgb]{0,0,0}\makebox(0,0)[lt]{\lineheight{1.25}\smash{\begin{tabular}[t]{l}$-\ln(1/\epsilon)$\end{tabular}}}}%
  \end{picture}%
\endgroup%

	\end{center}
\caption{The direction of the variable $t$ is thought of as the horizontal direction in the picture. The vertical directions correspond to the slicing of $M$ by 2-spheres.} 
	\label{fig.tubes1}
\end{figure}

Let us describe the sequence of Riemannian manifold we want to study. The manifold $M$ is a topological 3-sphere. The central part $M_C$ of $M$, in green in Figure \ref{fig.tubes1}, is the product $[-\ln(1/ \epsilon), \ln(1/ \epsilon)] \times \mathbb{S}^2$. Metrically, it is a warped product $g_{\epsilon} := dt^2 + \epsilon^2 \cdot \cosh(t)^2 \cdot g_{\mathbb{S}^2}$, where $g_{\mathbb{S}^2}$ denotes the round metric of curvature 1 on the 2-sphere. \\

The metric on some uniformly-in-$\epsilon$-fat-neighbourhood of $\partial M_C$, for example 
$$[-\ln(1\slash \epsilon), -\ln(1\slash \epsilon)+1] \cup [\ln(1\slash \epsilon)-1, \ln(1\slash \epsilon)] \times \mathbb{S}^2 \ , $$ 

 is smoothly close to the product metric on $\mathbb{S}^2 \times [0,1]$ and
interpolating with the round metric on a half $\mathbb S^3$, hence has uniformly bounded geometry with respect to $\epsilon$. Note also that the central part has Ricci curvature bounded from below but diameter $\ge -2 \ln(1\slash \epsilon)$. The metric does not need to be explicit on the extreme parts $M_L$ and $M_R$, in red in Figure \ref{fig.tubes1}: we only require that they have uniformly bounded geometry with respect to $\epsilon$. \\

As in the previous section, we shall emphasize with respect to which metric we measure by using the subscript $\epsilon$. In order to make the above described sequence of Riemannian metrics a counterexample, we show that $(\lambda^1_1)_{\epsilon}$ goes to $0$ as $\epsilon \to 0$ and that $(h_{\epsilon}^1)_{\epsilon \in ]0,1]}$ is bounded below by some constant independent of $\epsilon$. The following lemma addresses the eigenvalue part.

\begin{lemma}
	\label{lem.counterexample.eigenvalue}
	With the notation introduced above, 
	$$ (\lambda^1_1)_{\epsilon} \tends{\epsilon \to 0} 0 \ .$$
\end{lemma}

\textbf{Proof.} Since $H^1(\mathbb{S}^3, \RR) = \{0\}$, there is no harmonic $1$-forms and the conclusion of the above lemma follows from finding a coexact $1$-form with small Rayleigh quotient. We shall exhibit one supported on the central part of $M$. We denote by $g_{\mathbb{S}^2}$ the round metric on $\mathbb{S}^2$. Let $\alpha$ be any coclosed 1-form of $(\mathbb{S}^2, g_{\mathbb{S}^2})$ that we fix once and for all and $f$ be a function compactly supported on $] \ln(\epsilon), - \ln(\epsilon)[$. We set 
	$$ (\beta)_{(x,t)} := f(t) \alpha_x \ ,$$
which is a well defined 1-form supported on the central part of $M$. Since the metric is a warped product on the central part we have 
	$$ \star_{\epsilon} \alpha = (\star_{\mathbb{S}^2} \alpha) \wedge dt \ ,$$
and then  
	$$\star_{\epsilon} \beta = \ f(t) \  (\star_{\mathbb{S}^2} \alpha) \wedge dt $$ which shows that $\beta$ is coclosed since we assumed that $\alpha$ was coclosed with respect to $g_{\mathbb{S}^2}$. We then compute 
		$$ \beta \wedge \star_{\epsilon} \beta = f^2(t) \  \alpha \wedge (\star_{\mathbb{S}^2} \alpha) \wedge dt \ , $$
which after integration yields 
	$$ || \beta ||^2_{L^2_{\epsilon}} = || \alpha ||^2_{L^2(\mathbb{S}^2)} \inte{ \ln(\epsilon)}^{- \ln(\epsilon)} f^2(t) \ dt \ . $$ 
	
	Let us now compute the $L^2$ norm of $d \beta$. Since 
		$$d \beta = f'(t) \ dt \wedge \alpha + f(t) \ d \alpha $$ 
we have
		$$ \star_{\epsilon} d \beta = f'(t) \ dt \wedge \star_{\mathbb{S}^2} \alpha + \frac{f(t)}{ \epsilon^2 \cosh(t)^2} dt \ ,$$ 
and then 
	$$ d \beta \wedge \star_{\epsilon} d \beta = f'(t)^2 dt \wedge \alpha \wedge \star_{\mathbb{S}^2} \alpha + \frac{f^2(t)}{ \epsilon^2 \cosh(t)^2} dt \wedge d \alpha \ .$$ 

After integrating over $M$ and by noticing that the above second term vanishes (using Fubini's theorem and since $d \alpha$ is exact) we get
	$$ || d \beta ||^2_{L^2_{\epsilon}} = || \alpha ||^2_{L^2(\mathbb{S}^2)} \ \inte{ \ln(\epsilon)}^{- \ln(\epsilon)} f'(t)^2 \ dt \ .$$ 
In the end, we have shown 
$$ \frac{|| d \beta ||_{L^2_{\epsilon}} }{|| \beta ||_{L^2_{\epsilon}} } = \frac{|| f' ||_{L^2(] \ln(\epsilon), - \ln(\epsilon)[) } }{ || f ||_{L^2(] \ln(\epsilon), - \ln(\epsilon)[) }  } \ ,  $$
which corresponds to the energy of a compactly supported function $f$ on the interval $] \ln(\epsilon), - \ln(\epsilon)[$. One can then choose $f$ to realise the minimum: it corresponds to the first eigenfunction of the interval with Dirichlet boundary conditions. The associated eigenvalue is known to converge to $0$ as the length of the interval goes to infinity, which is the case as $\epsilon \to 0$.   \hfill $\blacksquare$ \\ 

The second part of the proof consists in showing that $h^1_{\epsilon}$ is uniformly bounded from below which will be an immediate consequence of 

\begin{proposition}
\label{prop.cunterexample.cheegerconstant}
There is a constant $C >0$ such that for any $\epsilon > 0$, for any embedded curve $\gamma$ in $\mathbb{S}^3$ there is a surface $S_{\gamma, \epsilon}$ such that 
	$$ \area_{\epsilon}(S_{\gamma, \epsilon}) \le C \cdot l_{\epsilon}(\gamma) \ . $$
\end{proposition}

The above proposition implies in particular that $ h^1_{\epsilon} \ge C^{-1}$ for all $\epsilon$, which yields to the desired conclusion. Its proof occupies the rest of this section. \\

\textbf{Proof.} In the sequel, $C$ will denote a constant independent of the curve $\gamma$ and of the parameter $\epsilon$ whose value may change from line to line. \\

We shall give more or less explicitly the surface $S_{\gamma, \epsilon}$. In order to do so, we will cut the curve $\gamma$ in suitable pieces $(\tilde{\gamma_i})_{i \in I}$ depending on when $\gamma$ enters the different parts of $M$. The endpoints of any of such pieces will then be glued back together to get a family of closed curves $(\gamma_i)_{i \in I}$ all of which belonging to some definite part, $M_L, M_C$ or $M_R$ of $M$. We shall see that there is a constant $C> 0$ such that for all $i \in I$ and all $\epsilon > 0$, there is a surface $S_{\gamma_i}$ whose boundary is $\gamma_i$ and which satisfies 
	$$ \area_{\epsilon}(S_{\gamma_i,\epsilon}) \le C \cdot l_{\epsilon}(\gamma_i) \ . $$
 We will then slit these surfaces accordingly on how we have cut $\gamma$ and paste them together to recover a surface $S_{\gamma, \epsilon}$ whose boundary is our initial curve $\gamma$ and which satisfies to the desired inequality. \\

\textbf{If the curve is contained in a definite part.} Let us first investigate the case where $\gamma$ will not need to be cut, namely the case where $\gamma$ is contained in some definite part. We start with the central part, which is the most important one since it is responsible for the large diameter of $(M, g_{\epsilon})$. \\

\begin{figure}[h!]
\begin{center}
	\def\svgwidth{0.8 \columnwidth}
%% Creator: Inkscape 1.0.1 (c497b03c, 2020-09-10), www.inkscape.org
%% Accompanies image file 'tubes2.pdf' (pdf, eps, ps)
\begingroup%
  \makeatletter%
  \providecommand\color[2][]{%
    \errmessage{(Inkscape) Color is used for the text in Inkscape, but the package 'color.sty' is not loaded}%
    \renewcommand\color[2][]{}%
  }%
  \providecommand\transparent[1]{%
    \errmessage{(Inkscape) Transparency is used (non-zero) for the text in Inkscape, but the package 'transparent.sty' is not loaded}%
    \renewcommand\transparent[1]{}%
  }%
  \providecommand\rotatebox[2]{#2}%
  \newcommand*\fsize{\dimexpr\f@size pt\relax}%
  \newcommand*\lineheight[1]{\fontsize{\fsize}{#1\fsize}\selectfont}%
  \ifx\svgwidth\undefined%
    \setlength{\unitlength}{1240.39039911bp}%
    \ifx\svgscale\undefined%
      \relax%
    \else%
      \setlength{\unitlength}{\unitlength * \real{\svgscale}}%
    \fi%
  \else%
    \setlength{\unitlength}{\svgwidth}%
  \fi%
  \global\let\svgwidth\undefined%
  \global\let\svgscale\undefined%
  \makeatother%
  \begin{picture}(1,0.38603765)%
    \lineheight{1}%
    \setlength\tabcolsep{0pt}%
    \put(0,0){\includegraphics[width=\unitlength,page=1]{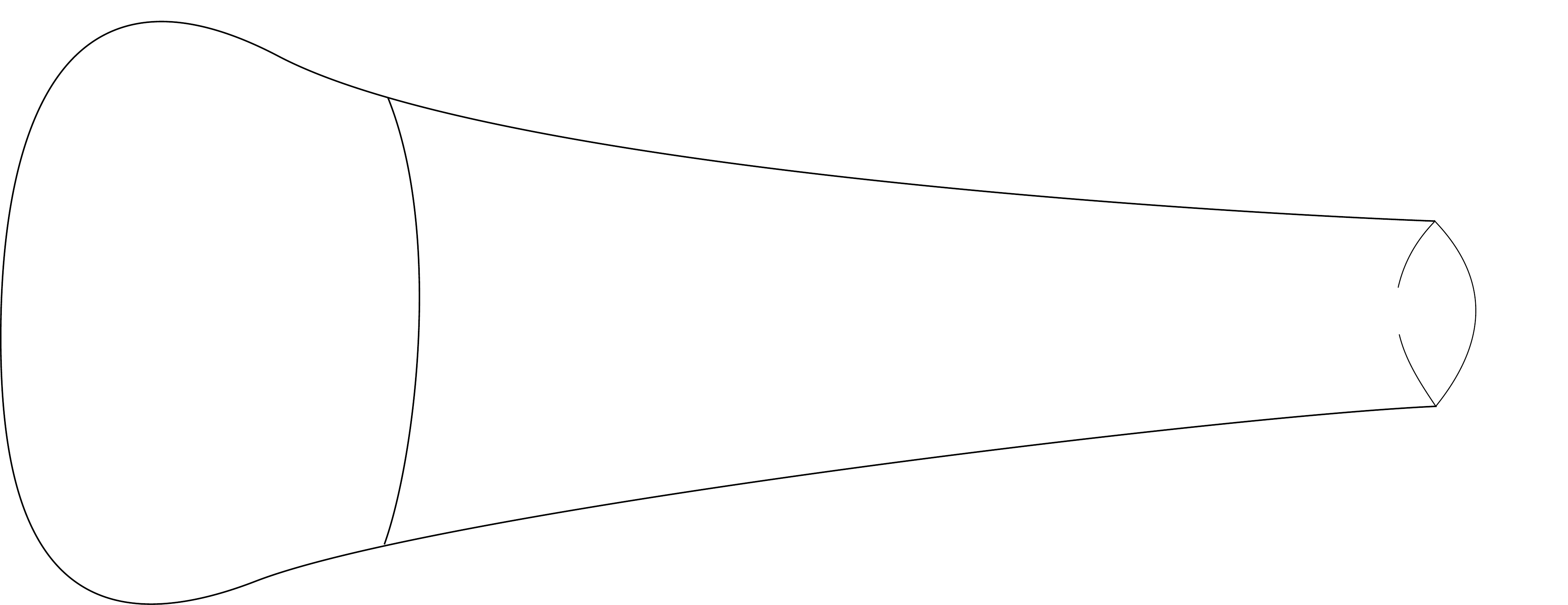}}%
    \put(0.88800084,0.09208303){\color[rgb]{0,0.36862745,0.14901961}\makebox(0,0)[lt]{\lineheight{1.25}\smash{\begin{tabular}[t]{l}$\mathbb{S}_0^2$\end{tabular}}}}%
    \put(0,0){\includegraphics[width=\unitlength,page=2]{tubes2.pdf}}%
    \put(0.51860046,0.36742586){\color[rgb]{0,0,1}\makebox(0,0)[lt]{\lineheight{1.25}\smash{\begin{tabular}[t]{l}$\mathbb{S}_{\gamma_0}$\end{tabular}}}}%
    \put(0.80151237,0.36571566){\color[rgb]{0,0,0}\makebox(0,0)[lt]{\lineheight{1.25}\smash{\begin{tabular}[t]{l}$\mathbb{\gamma}_0$\end{tabular}}}}%
    \put(0.34315762,0.16733204){\color[rgb]{0,0,0}\makebox(0,0)[lt]{\lineheight{1.25}\smash{\begin{tabular}[t]{l}$\gamma$\end{tabular}}}}%
    \put(0.56885639,0.17673817){\color[rgb]{0,0,0.8}\makebox(0,0)[lt]{\lineheight{1.25}\smash{\begin{tabular}[t]{l}$H(\gamma, \gamma_0)$\end{tabular}}}}%
    \put(0.08748211,0.16904224){\color[rgb]{0,0,0}\makebox(0,0)[lt]{\lineheight{1.25}\smash{\begin{tabular}[t]{l}$\huge{M_L}$\end{tabular}}}}%
  \end{picture}%
\endgroup%

	\end{center}
\caption{The 'whole homotopy' $H(\gamma, \gamma_0)$ is represented in light blue in this picture. The 'good surface' $S_{\gamma_0}$ is represented in dark blue. The seek surface $S_{\gamma}$ correspond to everything coloured in blue here.} 
	\label{fig.tubes2}
\end{figure}

We denote by $\gamma_0$ the projection of $\gamma$ on the 2-sphere $\mathbb{S}^2_{0}$ which lies right in the middle of the central part, see Figure \ref{fig.tubes2}. The idea to construct a surface whose boundary is $\gamma$ is to glue a 'good' $S_{\gamma_0}$ surface of $\mathbb{S}^2_{0}$ whose boundary is $\gamma_0$ together with 'the whole homotopy' $H(\gamma, \gamma_0)$ from $\gamma$ to $\gamma_0$, see Figure \ref{fig.tubes2}. \\

The surface $H(\gamma, \gamma_0)$ can be explicitly described by using the product structure of $M_C$: we parametrize $\gamma$ by
	$$ \gamma(s) := ( t(s), \gamma_0(s)) $$ where we denoted by $\gamma_0(s)$, slighly abusing the notation, the image of $\gamma$ under the projection $M = [\ln(\epsilon), -\ln(\epsilon)] \times \mathbb{S}^2 \to \mathbb{S}^2$ and where $s \in \mathbb{S}^1$ is the parameter. Then the 'whole homotopy' writes like
\begin{equation}
	\label{eq.wholehomotopy}
H(\gamma, \gamma_0) := \{ ( u \cdot t(s) ,   \gamma_0(s)) \ | \ \ u \in [0, 1] \ , \ s \in \mathbb{S}^1 \} \ .
\end{equation}

The 'good' surface is chosen as any surface on  $S_{\gamma_0} \subset \mathbb{S}^2 $ such that 
 $$ \area_{\mathbb{S}^2}(S_{\gamma_0}) \le C \cdot l_{\mathbb{S}^2}(\gamma_0) \  \  \text{ and} \  \ \partial S_{\gamma_0} = \gamma_0 $$
with $C$ which does not depend on $\gamma_0$. Note that we know such a constant $C > 0$ to exist for the same reason that we know that $h^1 > 0$. Note also that the metric on the $\mathbb{S}^2_0$ is $ \epsilon^2 \cdot g_\mathbb{S}^2$ with $\epsilon \le 1$, therefore:
 \begin{equation}
 	\label{eq.isopcounterexample}
 	  \area_{\epsilon}(S_{\gamma_0}) = \epsilon^2 \cdot \area_{\mathbb{S}^2}(S_{\gamma_0}) \le  \epsilon^2 \cdot C \cdot l_{\mathbb{S}^2}(\gamma_0) = \epsilon \cdot C \cdot l_{\epsilon}(\gamma_0) \le C \cdot  l_{\epsilon}(\gamma_0)  \ , 
 	\end{equation}
since $ \epsilon \le 1$. The point is that the constant appearing in the above upper bound can be chosen not to depend on $\epsilon$.  We then set 
	 $$S_{\gamma,\epsilon} := H(\gamma, \gamma_0) \cup S_{\gamma_0} \ , $$
which as a subset of $I \times \mathbb{S}^2$ does not depend on $\epsilon$.  The key property concerns $H(\gamma, \gamma_0)$ as testified in 

\begin{lemma}
\label{lem.counterexample.cheegerconstant}
	Let $\gamma$ be any curve in $M$ contained in its central part and let $\gamma_0$ be its central projection. There is a constant $C >0$ such that for all $\epsilon > 0$ we have 
		$$ \area_{\epsilon}(H(\gamma, \gamma_0)) \le C \cdot l_{\epsilon}(\gamma) \ .$$
\end{lemma}

The proof of the above lemma consists in using the explicit expression of the metric on the central part $M_C$ and the parametrization \eqref{eq.wholehomotopy} of $H(\gamma, \gamma_0)$. The details are left to the reader. \\

Because of Lemma \ref{lem.counterexample.cheegerconstant} and Inequality \eqref{eq.isopcounterexample}, the surface $S_{\gamma, \epsilon}$ verifies 
	$$ \area_{\epsilon}(S_{\gamma, \epsilon}) \le C \cdot l_{\epsilon}(\gamma) ,  $$
for some constant $C$ independent of $\epsilon$, concluding for central curves. \\

Let us now investigate the case of extremal curves. Let $\gamma$ be completely included in, say, the left part of $M$. Recall that this part is smoothly close to half a round sphere of radius $1$ independently of $\epsilon$. Let $\Phi$ be a diffeomorphism between $M_L$ and the round half $3$-sphere which preserves the metrics up to a constant, so that $\Phi(\gamma)$ is included in half of a round $3$-sphere. In order to fall back in the compact setting, we consider the double manifold, a full round 3-sphere of radius 1. As now embedded in the round 3-sphere, there is a surface $S_{\Phi(\gamma)}$ such that 	
	$$ \area_{\mathbb{S}^3}(S_{\Phi(\gamma)}) \le C \cdot l_{\mathbb{S}^3}(\Phi(\gamma)) \ .$$
		
\begin{figure}[h!]
\begin{center}
	\def\svgwidth{1 \columnwidth}
%% Creator: Inkscape 1.0.1 (c497b03c, 2020-09-10), www.inkscape.org
%% Accompanies image file 'tubes4.pdf' (pdf, eps, ps)
\begingroup%
  \makeatletter%
  \providecommand\color[2][]{%
    \errmessage{(Inkscape) Color is used for the text in Inkscape, but the package 'color.sty' is not loaded}%
    \renewcommand\color[2][]{}%
  }%
  \providecommand\transparent[1]{%
    \errmessage{(Inkscape) Transparency is used (non-zero) for the text in Inkscape, but the package 'transparent.sty' is not loaded}%
    \renewcommand\transparent[1]{}%
  }%
  \providecommand\rotatebox[2]{#2}%
  \newcommand*\fsize{\dimexpr\f@size pt\relax}%
  \newcommand*\lineheight[1]{\fontsize{\fsize}{#1\fsize}\selectfont}%
  \ifx\svgwidth\undefined%
    \setlength{\unitlength}{1473.17936088bp}%
    \ifx\svgscale\undefined%
      \relax%
    \else%
      \setlength{\unitlength}{\unitlength * \real{\svgscale}}%
    \fi%
  \else%
    \setlength{\unitlength}{\svgwidth}%
  \fi%
  \global\let\svgwidth\undefined%
  \global\let\svgscale\undefined%
  \makeatother%
  \begin{picture}(1,0.34079471)%
    \lineheight{1}%
    \setlength\tabcolsep{0pt}%
    \put(0,0){\includegraphics[width=\unitlength,page=1]{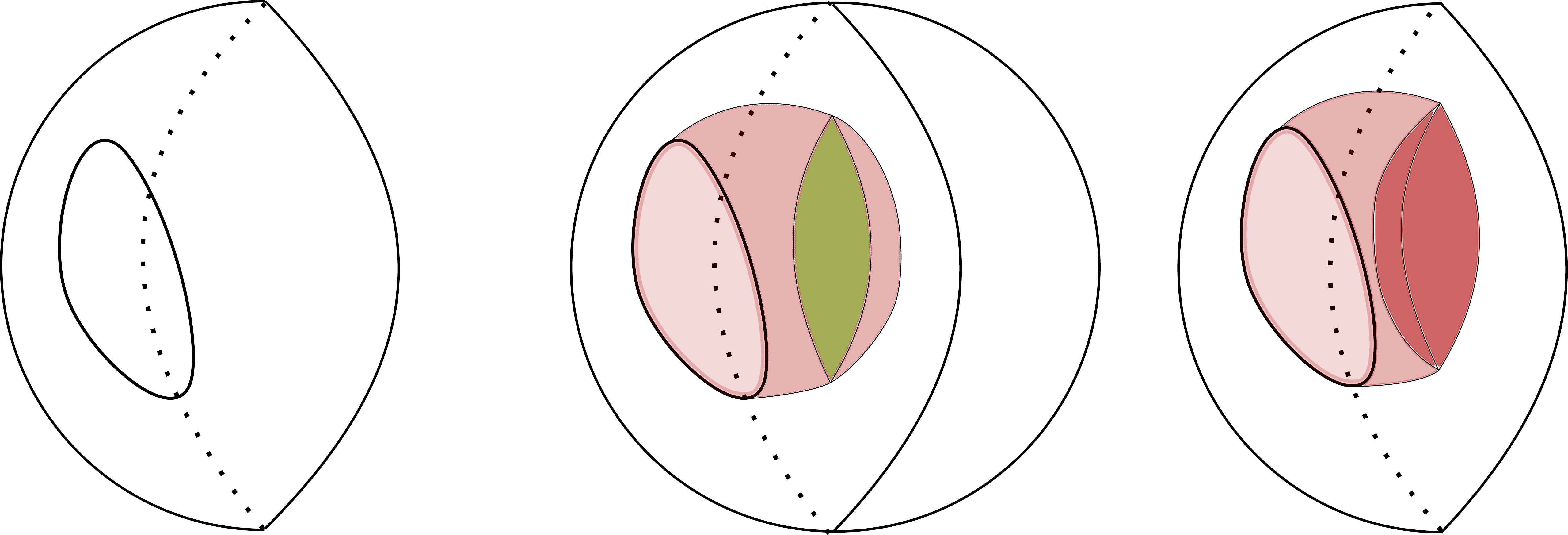}}%
    \put(0.11411387,0.18347868){\color[rgb]{0,0,0}\makebox(0,0)[lt]{\lineheight{1.25}\smash{\begin{tabular}[t]{l}$\Phi(\gamma)$\end{tabular}}}}%
    \put(0.79577641,0.06684187){\color[rgb]{0,0,0}\makebox(0,0)[lt]{\lineheight{1.25}\smash{\begin{tabular}[t]{l}$S_{\Phi(\gamma)}$\end{tabular}}}}%
    \put(0.42704543,0.0543763){\color[rgb]{0,0,0}\makebox(0,0)[lt]{\lineheight{1.25}\smash{\begin{tabular}[t]{l}$\widetilde{S_{\Phi(\gamma)}}$\end{tabular}}}}%
  \end{picture}%
\endgroup%

	\end{center}
\caption{On the left, the curve $\Phi(\gamma)$ in the half 3-sphere. On the middle, $\Phi(\gamma)$ seen as a curve of the full 3-sphere. As such there is a surface $\widetilde{S_{\Phi(\gamma)}}$ of controlled area whose boundary is $\Phi(\gamma)$. $\widetilde{S_{\Phi(\gamma)}}$ may not remain in the half part in which $\Phi(\gamma)$ is included. If not, we use the symmetry defining the doubling to get a surface $S_{\Phi(\gamma)}$, in red on the rightmost picture, of same area but included in the same half sphere than $\Phi(\gamma)$.} 
	\label{fig.tubes4}
\end{figure}
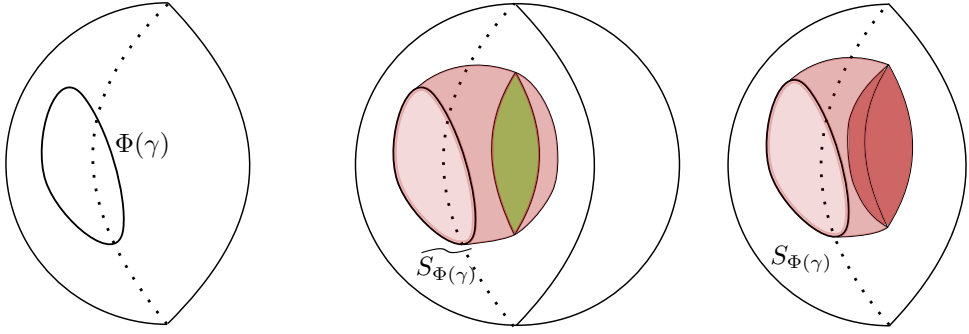

Using if necessary the symmetry between the two half 3-sphere, see Figure \ref{fig.tubes4}, one can assume that $S_{\Phi(\gamma)}$ is completely included in our original half sphere. Since $\Phi$ is an quasi-isometry, we have 
	$$ \area_{\epsilon}( \Phi^{-1}(S_{\Phi(\gamma)}))  \le C' \cdot \area_{\epsilon}(S_{\Phi(\gamma)}) \ \text{ and } \ l_{\epsilon}(\Phi^{-1}(\gamma)) \le C'' \cdot l_{\epsilon}(\gamma) \ . $$
In other word, 
\begin{equation}
	\label{eqdemcounterexample1000} 
 S_{\gamma} := \Phi^{-1}(S_{\Phi(\gamma)})
 \end{equation}
satisfies to the conclusion of Proposition \ref{prop.cunterexample.cheegerconstant}. \\

\textbf{Cutting off $\gamma$ and pasting altogether the $
S_{\gamma_i}$.} It remains then to do the cut and paste. As already mentioned, we shall cut the curve $\gamma$ whenever $\gamma$ goes from the central part to the extremal ones. Namely, the two 2-spheres corresponding to the boundary of the central part $\partial M_C$ will be considered as 'cutting subsurfaces': we cut $\gamma$ whenever it crosses some definite small neighbourhood of these. We are left with a collection of curves whose connected components $(\widetilde{\gamma_i})_{i \in I}$ are all contained in one of the 3 parts $M_L, M_C, M_R$ of $M$. This connected component will not be closed, expect in the already-dealt-with-case where $\gamma$ was initially contained in a definite part of $M$; we close them by considering any minimizing geodesic that relates the exit points of $\gamma$ to the re-entering points in $M_C$.\\

\begin{figure}[h!]
\begin{center}
	\def\svgwidth{1 \columnwidth}
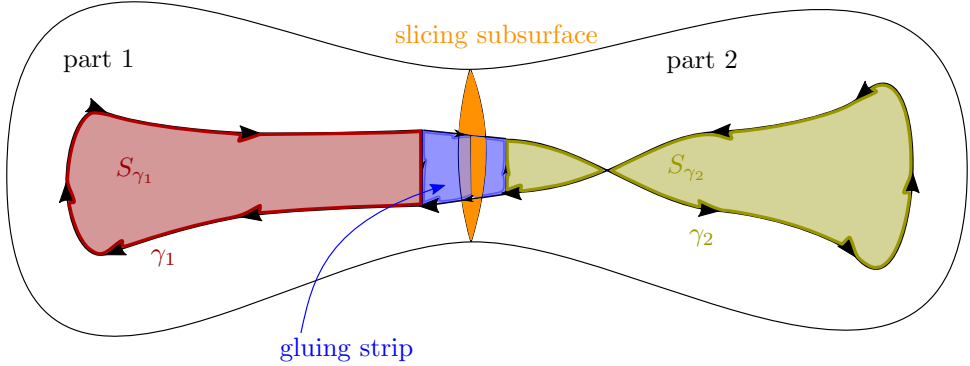
	\end{center}
\caption{In this schematic example, there is only one cutting subsurface, in orange, which cuts the manifold in two pieces, part 1 and part 2 (on the contrary of $M$ which has three parts). The black curve plays the role of $\gamma$. The red and green curves $\gamma_1$ and$\gamma_2$ correspond to the leftovers curves obtained after the cutting, closed up by short geodesics segments. These 2 closed curves are the boundaries of two surfaces $S_{\gamma_1}, S_{\gamma_2}$ that we glue back together with the help of a gluing strip in blue.} 
	\label{fig.tubes8}
\end{figure}

Recall that any definite neighbourhood of the cutting subsurfaces is of bounded geometry (independently of $\epsilon$). In particular, since we have chosen to close the curves with minimizing geodesics there is a constant $C > 0$, which does not depend on $\epsilon$, such that 
\begin{equation}
	\label{eqdemcounterexample1000} 
\somme{i \in I} \ l_{\epsilon}( \gamma_i) \le C \cdot l_{\epsilon}(\gamma) 
\end{equation}

Because the curves $(\gamma_i)_{i \in I}$ belong to some definite part of $M$, we can use the work already performed and find for any $i \in I$ a surface $S_i$ such that $ \area_{\epsilon}(S_i) \le C l_{\epsilon}(\gamma_i) $ for some constant $C > 0$ independent of $\epsilon$. We then glue back all these subsurfaces together according to how we cut the curve $\gamma$ by using 'gluing strips' (see Figure \ref{fig.tubes2}). The resulting surface has, by construction, $\gamma$ as boundary. Let us conclude by seeing that its area satisfies to the desired inequality of Proposition \ref{prop.cunterexample.cheegerconstant}. We have by construction
\begin{equation}
	\label{eqdemcounterexample10000}
	 \area_{\epsilon}(S_{\gamma}) = \somme{i \in I} \area_{\epsilon}(S_{\gamma_i}) + \area_{\epsilon}( \text{strips} ) \ .
\end{equation}

The area corresponding to the strips may be taken as small as desired, depending on how thin we choose them, and the corresponding terms can be ignored. Since we know that $ \area_{\epsilon}(S_i) \le C \cdot  l_{\epsilon}(\gamma_i) $, we get for all $\epsilon$ 
\begin{equation}
	\label{eqdemcounterexample100000}
	 \area_{\epsilon}(S_{\gamma}) \le C \cdot  \somme{i \in I} l_{\epsilon}(\gamma_i) \ .
\end{equation}
 
Inequality \eqref{eqdemcounterexample1000} then gives 
\begin{equation*}
	\label{eqdemcounterexample100000}
	 \area_{\epsilon}(S_{\gamma}) \le C \cdot l_{\epsilon}(\gamma) \ , 
\end{equation*}
concluding. \hfill $\blacksquare$

\section{Concluding remarks and questions}
\label{sec.remarks}

1. One may wonder if a Cheeger-like inequality as the one given by the conclusion of Theorem \ref{maintheo} also holds for differential forms in any degree $p$. Indeed, the definition of $h^1$ may straightforwardly generalised to any degree: given a homologically trivial $p$-cycle we define its 'area' as 
	$$ \mu^{p+1}({d^{-1} S}) := \infi{ \partial N = S} \ \mu^{p+1}( N) $$ 
	and the associated isoperimetrical constant as 
	$$ h^p := \infi{ S \hookrightarrow M ,\  [S] = 0 } \ \frac{\mu^p(S)}{ \mu^{p+1}({d^{-1} S})}  \ .  $$
The main difficulty to overcome, as it was the case for 1-forms, is to find yet another analogous to the level sets. The first not already treated case concerns coexact 2-forms in dimension 5. \\ 

2. We discuss here, following up with the remark (2) stated after our main theorem \ref{maintheo}, the definition of another Cheeger's constant for 1-forms which may come with a statement analogous to \ref{maintheo} but with a constant $C$ which does not depend on the geometry of $M$, as in Cheeger's theorem. In a forthcoming article we will state 
such a Cheeger inequality between the following constant and $\lambda^1_1$ when $M$ has dimension 3:
	$$ h_*^1 := \infi{ \gamma \hookrightarrow M ,\  [\gamma] =0 } \frac{l(\gamma)}{|| \omega_{\gamma}||_{1}} \ ,  $$
where the infinum runs over all real homologically trivial smooth closed curves $\gamma$ and where $l(\gamma)$ is the length of $\gamma$ and where $\omega_{\gamma}$ is the magnetic linking form associated to $\gamma$. This 1-form, as suggested by its name, has a strong topological flavour, when considered as the Hodge dual of a meridian of $\gamma$, and a physical one, when interpreted as the magnetic field generated by a stationary electrical current flowing through $\gamma$. This physical/topological relation was discovered by Gauss in the 19th century. 
%We will see that a Cheeger like inequality holds with the above defined constant in dimension 3 without any dependency on the geometry of $M$. 
%The definition of $h_*^1$ may seem artificial since the denumerator $|| \omega_{\gamma}||_{1}$ does not seem directly related to the geometry of $M$. 
In the definition of $h_* ^1$, the denominator can actually be interpreted as the average area of the pages of an 'open book decomposition' of $M \setminus \gamma \to \gamma$, so that
the minimal area of a surface with boundary $\gamma$ in the definition of $h^1$ is replaced here by an average.
However, the constant $h_*^1$ is more analytical by nature since defining $\omega_{\gamma}$ involves the Green operator. \\

3. It would be interesting to find examples where the quantity $h^1$ may be computed explicitly. The study performed in Section \ref{sec.example} suggests that, in the case of the round 3-sphere, the constant $h^1$ is realized by the isoperimetrical ratio of any of the great circle of $\mathbb{S}^3$.

\appendix

\section{positivity of $h_1$}
The goal of this section is to show that the constant $h_1$ is positive for any Riemannian manifold through the following

\begin{proposition}
	\label{prop_h1positive}
		Let $M$ be a closed Riemannian manifold. There is a constant $C >0$ such that for any smooth and real homologically trivial curve $\gamma$ of $M$ there is a rectifiable surface $S_{\gamma}$ such that
			\begin{equation}
				\label{eq.proph1positive}
					 \area(S_{\gamma}) \le C \cdot l(\gamma) \ .
			\end{equation}
\end{proposition} 

The above proposition implies in particular that $A(\gamma) \le C \cdot l(\gamma)$ by definition of $A(\gamma)$. We decided to state Proposition \ref{prop_h1positive} this way to emphasize that we will 'exhibit' a surface of controlled area in order to deduce an upper bound of $A(\gamma)$ rather than trying to compute it directly. \\

\textbf{Proof.} As in Section \ref{sec.counterexample}, the value of the constant $C$ appearing all along this section may vary from line to line. \\

 The proof relies on the existence of triangulations for smooth manifolds. We fix from now a triangulation of $M$ for which we denote by $\tau$ the associated simplicial complex and by $\tau_i$ its $i$-skeleton. The first step of the proof consists in reducing the study to the combinatorial case of $\gamma$ belonging to the 1-skeleton of $\tau$. 

\begin{lemma}
	\label{lemma.reduc1skeleton}
		%Under the assumption of Proposition \ref{prop_h1positive}. 
		Let $M$ be a $d$-dimensional close manifold with a triangulation $\tau$.
		There is a constant $C > 0$ such that for any curve $\gamma$ of $M$ there is a curve $\gamma_1 \in \tau_1$ and a surface $S_{\gamma \to \gamma_1}$ such that
		\begin{equation}
			\left\{  
				\begin{array}{l}	
					\partial S_{\gamma \to \gamma_1} = \gamma \cup - \gamma_1  \\
					\area(S_{\gamma \to \gamma_1}) \le C \cdot l(\gamma) \\		
					l(\gamma_1) \le C \cdot l(\gamma) \ ,
				\end{array}
			\right.
		\end{equation}
where $ -\gamma_1$ stands for the curve $\gamma_1$ endowed with the opposite orientation.
\end{lemma}
The second step of the proof will consist in dealing with the remaining case where $\gamma$ belongs to $\tau_1$. 

\begin{lemma}
\label{lemma.1skeleton}
	%Under the assumption of Proposition \ref{prop_h1positive}. 
	%Let $M$ be a $d$-dimensional close manifold with a triangulation $\tau$.
	Let $M$ be a $d$-dimensional close manifold with a triangulation $\tau$. There is a constant $C > 0$ such that the following holds. For any curve real homologically trivial $\gamma_1 \subset \tau_1$, that is, such that 
	$r \cdot \gamma _1 =0$ in $H_1 (M, \mathbb Z)$,
	there is a surface $S_{\gamma_1}$ satisfying
		$$	\area(S_{\gamma_1}) \le C \cdot l(\gamma_1) \ \ \text{ and } \ \ \partial S_{\gamma_1} = r \cdot\gamma_1.$$
\end{lemma}

For a given curve $\gamma$, it is easy to verify that the surface 
	$$ S_{\gamma} :=  S_{\gamma \to \gamma_1} \cup S_{\gamma_1} \ , $$
where $\gamma_1, S_{\gamma \to \gamma_1}$ and $S_{\gamma_1}$ are as in the conclusions of  Lemmas  \ref{lemma.reduc1skeleton}  and \ref{lemma.1skeleton},
 satisfy Inequality \eqref{eq.proph1positive}. One is then left to prove the two above mentioned lemmas. \\

\textbf{Proof of Lemma \ref{lemma.reduc1skeleton}.} We will have to work with a slightly greater class of metric spaces than the one of manifolds, which is why we restate Lemma \ref{lemma.reduc1skeleton} as follows. \\

\textit{Let $\tau$ be a compact metric simplicial complex of dimension $d$, that is, endowed with a continuous Riemannian metric which is smooth in restriction to the interiors of 
the $i$-skeletons. There is a constant $C > 0$ such that for any curve $\gamma$ of $M$ there is a curve $\gamma_1 \in \tau_1$ and a surface $S_{\gamma \to \gamma_1}$ such that}
		\begin{equation}\label{recur}
			\left\{  
				\begin{array}{l}	
					\partial S_{\gamma \to \gamma_1} = \gamma \cup -\gamma_1  \\
					\area(S_{\gamma \to \gamma_1}) \le C \cdot l(\gamma) \\		
					l(\gamma_1) \le  C \cdot l(\gamma) \ .
				\end{array}
			\right. 
		\end{equation}

Note that in particular the assumptions above encompass the case of a triangulable Riemannian manifold. The proof of the above statement goes by induction on the dimension $d$. Let then $\tau$ be a metric simplicial complex as above and $\gamma$ be a curve of $\tau$. If the curve $\gamma$ is in the $1$-skeleton $\tau _1$, there is nothing to prove. We thus assume that (\ref{recur}) is proved for any curve in the $(d-1)$-skeleton $\tau _{d-1}$ and  let us prove it
for $\gamma$ in $\tau _d$.\\

We split $\gamma$ accordingly to when it enters/exits two different $d$-dimensional simplex of $\tau_d$. We are then left with a collection of possibly non closed curves $(\gamma^j)_{j \in J}$ whose endpoints lie in $\tau_{d -1}$. All the simplex $(\tau_d^i)_{i \in I}$ of $\tau_d$ are diffeomorphic to the euclidean simplex $ \Delta_d$ of dimension $d$. Since $M$ is compact, there is only finitely many $d$-simplex and therefore they are uniformly Lipschitz equivalent to $\Delta_d$. We rely on the euclidean local analogous of Lemma \ref{lemma.reduc1skeleton}.

 \begin{lemma}
 	\label{lemma.simplex.curve}
 	Let $\Delta_d$ be the Euclidean simplex of dimension $d$. There is a constant $C$ which depends only on $d$ such that for any curve $\gamma$ of $\Delta_{d}$ whose endpoints are in $\partial \Delta_d$ there is a curve $\gamma_{d-1} \in \partial \Delta_{d}$ and a surface $S_{\gamma \to \gamma_{d-1}}$ such that
		\begin{equation}
			\left\{  
				\begin{array}{l}	
					\partial S_{\gamma \to \gamma_{d-1}} = \gamma \cup -\gamma_{d-1}  \\
					\area(S_{\gamma \to \gamma_{d-1}}) \le C \cdot l(\gamma) \\		
					l(\gamma_{d-1}) \le C \cdot  l(\gamma) \ .
				\end{array}
			\right.
		\end{equation}
 \end{lemma}
 
The proof of the above lemma is classical: one has to choose a minimizing geodesic $\gamma_{d-1}$ of $\partial \Delta_d$ whose endpoints are the same as those of $\gamma$. We are then left with a closed curve $\gamma \cup - \gamma_{d-1}$ for which we know that there is a surface $S_{\gamma \to \gamma_{d-1}}$ such that
$$\area(S_{\gamma \to \gamma_{d-1}}) \le C \cdot (l(\gamma) + l(\gamma_{d-1})) \ . $$
Such a surface can be constructed explicitly by fixing a base point in $\Delta$ and by relating all the point of the closed curve $\gamma \cup - \gamma_{d-1}$ to it by euclidean segments. This construction (or variant of it) is sometimes referred as the 'cone construction', see for example \cite{artfedererflemming}. We conclude by noticing that $l(\gamma_{d-1}) \le C \cdot l(\gamma)$ for some constant $C$ since $\gamma_{d-1}$ was taken as a minimizing geodesic of $\tau_{d-1}$. The detail of the proof are left to the reader. \\ 

Let us see how to use the above lemma to fall back within the induction assumption. Let us denote by $(\tilde{\gamma_i})_{i \in I}$ the images in $\Delta_d$ under the diffeomorphisms between $(\tau^i_d)_{i \in I}$ and $\Delta_d$ of the curves $(\gamma_i)_{i \in I}$. We use to above lemma with the $(\tilde{\gamma})^j_{j \in J}$. Because all the $\tau^i_d$ are uniformly quasi-isometric to $\Delta _d$ there is a constant $C > 0$ such that for any $j \in J$ there is a curve $\gamma^j_{d-1}$ and a surface $S_{\gamma^j \to \gamma^j_{d-1}}$ such that 
\begin{equation}
			\left\{  
				\begin{array}{l}	
					\partial S_{\gamma^j \to \gamma^j_{d-1}} = \gamma^j \cup - \,\gamma^j_{d-1}  \\
					\area(S_{\gamma^j \to \gamma^j_{d-1}}) \le C \cdot l(\gamma^j) \\		
					l(\gamma^j_{d-1}) \le C \cdot l(\gamma^j) \ .
				\end{array}
			\right.
		\end{equation}

We set
	$$ \gamma_{d-1} := \underset{j \in J}{\cup} \gamma_{d-1}^j \ \ \text{ 	and  } \  \ S_{ \gamma \to \gamma_{d-1}} :=  \underset{j \in J}{\cup} S_{\gamma^j \to \gamma^j_{d-1}} \ . $$

Because of lemma \ref{lemma.simplex.curve}, both the length of the curve $\gamma_{d-1}$  and the area of $S_{ \gamma \to \gamma_{d-1}}$ are controlled by the length of $\gamma$. We now use the induction assumption with $\gamma_{d-1}$ to get a a curve $\gamma_1$ and a surface $S_{\gamma_{d-1} \to \gamma_1}$ which satisfies all the conclusion of Proposition \ref{prop_h1positive}. We conclude by setting 
	$$ S_{\gamma \to \gamma_1} := S_{ \gamma \to \gamma_{d-1}}  \cup S_{\gamma_{d-1} \to \gamma_1} \ , $$
	which, together with the curve $\gamma_1$ given by the induction assumption, satisfies all the requirements of Proposition \ref{prop_h1positive}. \hfill $\blacksquare$ \\
	
	Let us conclude this section by proving Lemma \ref{lemma.1skeleton}. \\

\textbf{Proof of Lemma \ref{lemma.1skeleton}.} Because $M$ is compact, the triangulation $\tau$ is finite. In particular, the vectorial space given by the 1-chains of $\tau$ is finite dimensional and comes with a natural basis $e_1, ... , e_p$ given by the edges of the triangulation. \\

Since the boundary operators $\partial$ are linear (in particular the one defined over the 2 chains into the 1-chains), the vectorial space of homologically trivial 1-chains $\partial \tau_2 \subset \tau_1$ is also finite dimensional. We fix once and for all a basis $\nu_1, ... , \nu_k$ of $\partial \tau_2$ that we complete with $\nu_{k + 1} , ... , \nu_p$ as a basis of $\tau_1$. \\

Note that there is a constant $C > 0$ such that the coefficients of 1-chain $\gamma_0$ expressed in the basis of $\tau_1$ given by the edges writes like 
	$$ \gamma_0 = \somme{1 \le i \le p} a_i e_i \ , $$
then 
$$ \somme{1 \le i \le p} |a_i| \le C \cdot l(\gamma_0)   \ . $$
Since $\tau_1$ is finite, up to the choice of a linear change of basis, one has as well that if a 1-chain $\gamma$ expressed in the basis given by the $\nu_i$ writes like 
	$$ \gamma = \somme{1 \le i \le p} b_i \nu_i \ , $$
then 
\begin{equation}	
	\label{eq.appendix1}
\somme{1 \le i \le p} |b_i| \le  C \cdot  l(\gamma) \ , 
\end{equation}
for some constant $C > 0$ (the $L_1$ norms associated to the basis $(e_i)$ and $(\nu_i)$ are equivalent). \\

For any $1 \le i \le k$ we now fix a 2-chain $S_i$ such that $\partial S_i = \nu_i$. Note in particular that we have for any $1 \le i \le k$ 
	$$ A( \nu_i) \le \supr{1 \le i \le k} \area(S_i) \ .$$

If $\gamma_0 \in \partial \tau_2$ it decomposes with respect to the $(\nu_i)_{1 \le i \le k}$:
	$$ \gamma_0 = \somme{1 \le i \le k} b_i \nu_i \ . $$

In particular, the 2-chains 
	$$ S_{\gamma_0} := \underset{1 \le i \le k}{\cup} b_i S_i $$
has $\gamma_0$ has boundary. By construction we have 
	\begin{align*}
		 A(\gamma_0) \le 	\area(S_{\gamma_0}) = \somme{1 \le i \le k}  |b_i| \area(S_i)  \le C \cdot \somme{1 \le i \le k}  |b_i| \ , 
	\end{align*}
where $C := 	 \supr{1 \le i \le k} \area(S_i)$. Combined with Inequality \eqref{eq.appendix1} we get
	$$  A(\gamma_0) \le C \cdot l(\gamma_0) \ ,$$ 	concluding.   \hfill $\blacksquare$ $\blacksquare$

\bibliographystyle{alpha}

\bibliography{bibliography} 

\bigskip
\textsc{A. Boulanger, Dipartimento di Matematica, Università di Bologna, Via Zamboni 33, 40126 Bologna, Italy}\par\nopagebreak
  \textit{E-mail address}: \texttt{adrien.boulanger@unibo.it}

\medskip

\textsc{G. Courtois, Institut de Math\'ematiques de Jussieu-Paris Rive Gauche, CNRS et Universit\'e Pierre et Marie Curie, 4 place Jussieu, 75232 Paris Cedex 09, France}\par\nopagebreak
  \textit{E-mail address}: \texttt{gilles.courtois@imj-prg.fr}

\end{document}